\documentclass{elsarticle}

\usepackage{setspace}
\usepackage{amsmath}
\usepackage{amssymb}
\usepackage{amsfonts}
\usepackage{mathrsfs}
\usepackage{bm}
\usepackage{bbm}
\usepackage[title]{appendix}
\usepackage{graphicx}
\usepackage{tikz}
\usepackage{pgfplots}
    \usetikzlibrary{spy}
\usepackage{caption}
\usepackage{subcaption}
\pgfplotsset{compat = 1.3}
\usepackage[margin=1in]{geometry}
\usepackage{comment}
\usepackage{multirow}
\usepackage{xcolor}

\definecolor{myyellow}{rgb}{0.75, 0.56, 0.01}
\definecolor{mygreen}{rgb}{0.21, 0.45, 0.11}

\journal{Journal of Computational Physics}

\begin{document}

\begin{frontmatter}

\title{Half-closed discontinuous Galerkin discretisations}

\author[rvt2,rvt3]{Y.~Pan\fnref{fn1}\corref{cor1}}
\ead{yllpan@berkeley.edu}

\author[rvt2,rvt3]{P.-O.~Persson\fnref{fn3}}
\ead{persson@berkeley.edu}

\address[rvt2]{Department of Mathematics, University of California, Berkeley, Berkeley, CA 94720, United States}
\address[rvt3]{Mathematics Group, Lawrence Berkeley National Laboratory, 1 Cyclotron Road, Berkeley, CA 94720, United States}
\cortext[cor1]{Corresponding author}
\fntext[fn1]{Graduate student, Department of Mathematics, University of California, Berkeley}
\fntext[fn3]{Professor, Department of Mathematics, University of California, Berkeley}


\begin{abstract}
We introduce the concept of half-closed nodes for nodal discontinuous Galerkin (DG) discretisations. Unlike more commonly used closed nodes in DG, where on every element nodes are placed on all of its boundaries, half-closed nodes only require nodes to be placed on a subset of the element's boundaries. The effect of using different nodes on DG operator sparsity is studied and we find in particular for there to be no difference in the sparsity pattern of the Laplace operator whether closed or half-closed nodes are used. On quadrilateral/hexahedral elements we use the Gauss-Radau points as the half-closed nodes of choice, which we demonstrate is able to speed up DG operator assembly in addition to leverage previously known superconvergence results. We also discuss in this work some linear solver techniques commonly used for Finite Element or discontinuous Galerkin methods such as static condensation and block-based methods, and how they can be applied to half-closed DG discretisations.
\end{abstract}

\end{frontmatter}
	
\section{Introduction}
The discontinuous Galerkin (DG) method is a popular variant of the Finite Element method which allows for discontinuities between elements in the solution space \cite{cockburn2001rkdg,hesthaven2007nodal}. Some properties of the method include its ability to easily achieve arbitrarily high orders of accuracy for a wide range of problems on unstructured meshes and to naturally allow for stabilisation through the use of approximate Riemann solvers. This has made it an attractive proposition for many fluid dynamics applications and it has as a result garnered much research attention in the past few decades since its introduction by Reed and Hill \cite{reed1973triangular} for the neutron transport problem.

One common criticism of the DG method however is its cost. In addition to having extra degrees of freedom compared to continuous Finite Elements, operator assembly in DG is often more expensive as numerical quadrature is necessary both inside element volumes and on their boundaries. This has led to significant research over the years on the development of related methods to address this, such as the DG-Spectral Element Method \cite{kopriva2010quadrature,kopriva1996dgsem},  line-based DG \cite{linedg}, Spectral Volumes \cite{wang2002sv}, Spectral Differences \cite{liu2006sd}, Flux Reconstruction \cite{huynh2007fr}, and Hybridized Discontinious Galerkin \cite{cockburn2009unified}. Many of these methods either alter the DG formulation or consider the equations in differential form instead of an integral one, so as to bypass some of the costs associated with the method.

In this work, we do not change the underlying DG formulation but instead introduce the use of half-closed nodes in DG discretisations to address this problem of cost. These nodes differ from commonly used closed nodes such as the Gauss-Lobatto or Chebyshev nodes of the second kind in DG and its related methods, where within each element the nodes are placed on all its boundaries creating the trademark node duplication pattern on inter-element boundaries. Closed nodes are often chosen to ensure direct coupling between neighbouring elements and to allow for efficient computation of boundary terms. Half-closed nodes on the other hand only require that on each inter-element boundary between any two elements, only one of the elements must have nodes placed on the boundary. This allows for increased flexibility in node placement and we show that this simple modification can enable more efficient assembly of DG operators without modifying the underlying DG formulation and its properties.

For quadrilateral and hexahedral elements in particular we use the Gauss-Radau nodes as the half-closed nodes of choice. These nodes double as quadrature nodes for numerical integration similar to the Gauss-Lobatto nodes in the DG-Spectral Element Method, but attain an extra order of quadrature precision over Gauss-Lobatto nodes \cite{Abramowitz1972handbook}. This produces a diagonal mass matrix which we show can be used to efficiently assemble other DG operators such as the divergence and Laplace operators. To determine which subset of the element boundaries half-closed nodes are placed on, we introduce a simple procedure based on switch functions.

To get a clearer picture of the cost associated with using half-closed nodes, we analyse the effect of using half-closed nodes on the sparsity pattern of DG operators for both first and second order partial differential equations. For second order equations the Local Discontinuous Galerkin (LDG) method is employed. While for first order equations a small increase in the number of non-zero entries is observed in first derivative DG operators, the sparsity pattern for second derivative operators is identical whether using closed or half-closed nodes. This suggests the possibility for commonly used linear solver techniques for DG methods with closed nodes to be applicable also to DG methods with half-closed nodes.

To this effect we discuss two types of linear solvers. The first is that of static condensation, also known as Guyan reduction \cite{guyan1965reduction}, popular in Finite Element methods \cite{weaver1987structural}. The second is that of block-based techniques such as block-Jacobi, popular for DG methods \cite{persson2008newtongmres}. We consider the use of these two methods in this setting and show both these techniques to be readily applicable to both first and second derivative DG operators using either closed and half-closed nodes.

We finally perform a series of numerical tests to demonstrate the half-closed DG method as well as the application of the described linear solvers to the resulting operators. To test the half-closed DG method we apply the method to solving the advection and Poisson equations in 1D and 2D on quadrilateral meshes. For these equations we verify that the solution obtained at the half-closed Gauss-Radau nodes display the superconvergence properties established in \cite{adjerid2002posteriori,liu2021superconvergence, yang2015analysis}. To test the linear solvers we apply static condensation and block based methods to the LDG Laplace operator, and study their effects on its spectrum.

The paper is structured as follows. We first give a brief review of the DG method and the LDG framework for second order equations in Section 2. In Section 3, the sparsity patterns of both first and second order DG operators are examined for the different types of nodes, and the cost of assembling each of these operators analysed in Section 4. Specific details of half-closed nodes and their placement including the switch function construction are given in Section 5. In Section 6 we pivot to the topic of the linear solvers before performing some numerical examples in Section 7 and finally concluding in Section 8.

\section{Discontinuous Galerkin review}
We briefly review the details of the Discontinuous Galerkin (DG) method. We focus here on the method applied to first and second order equations, where for second order equations we 
use the Local Discontinuous Galerkin (LDG) method \cite{ldg1998}. Some resources for a more in depth discussion of the DG method include \cite{cockburn2012discontinuous, hesthaven2007nodal}.

\subsection{Discretisation} \label{sect:discretisation}
To formulate the DG method, a domain $\Omega$ is first discretised into elements $\mathcal{T}_h = \{ K_n : \bigcup_n K_n = \Omega\}$. A finite element space $V_h$ on these elements is introduced as
\begin{equation}
    V_h = \{ v_h : v_h|_{K_n} \in V(K_n), ~\forall n \}
\end{equation}
where $V(K_n)$ is a function space isolated to element $K_n$. A popular choice for $V(K_n)$ are polynomial spaces of degree at most $p$ which we focus on solely for this paper. For simplex elements they are commonly chosen to be the space of multivariate polynomials and for quadrilateral elements as the $d$-dimensional direct product of 1-dimensional polynomials.

Functions in $V_h$ are in general discontinuous along the boundaries of elements $\partial K_n$ as functions in $V(K_n)$ have local support confined to a single element. Consequently numerical flux functions need to be defined on element boundaries in order to properly define the value of a function in $V_h$ on these boundaries. Throughout the text, these numerical fluxes are denoted using hats $\hat{F}$, wherein the exact form of the numerical flux in general depends on the equation under consideration.

To represent the solution $u \in V_h$, a representative basis function set for $V(K_n)$ is introduced on each element, denoted $\{ \phi_j \}$. A common way to do this is to on each element introduce a set of nodes $s_j$, and to define an interpolatory set of functions such that $\phi_i(s_j) = \delta_{ij}$. For polynomial function spaces, the number of nodes needed in each element $N$ to define a full set of basis functions depends on the desired polynomial degree of functions in each element, and on the shape of the element itself. Specifically for a chosen polynomial degree $p$ in $d$-dimensions, the number of nodes required in a simplex element is $N={p+d \choose d}$, whilst in quadrilateral elements is $N=(p+1)^d$.

\subsection{Nodes} \label{sect:nodes}
Nodes are commonly placed in the DG method such that $N_b$ nodes lie on every boundary of the element so that they define a $(d-1)$-dimensional complete polynomial space of the same degree $p$ on the boundary. For polynomial function spaces in simplex elements this means that $N_b = {p+d-1 \choose d-1}$ nodes are placed on each element boundary, whilst in quadrilateral elements $N_b=(p+1)^{d-1}$ nodes are placed on each element boundary. Nodes placed in this manner are known as closed nodes. Closed nodes are chosen usually to minimise communication across elements, as functions on element boundaries can be evaluated through only the $N_b$ nodal values on that boundary, rather than requiring all nodal values of the element.

A less popular choice in DG is to use open nodes, where fewer than $N_b$ nodes are placed on all element boundaries. This extra flexibility in node placement enables the use of high precision quadrature points such as the Gauss-Legendre points as solution nodes. The tradeoff however is that when using open nodes increased density in the off-diagonal blocks in the DG operators is generally observed, as evaluation of values on an element boundary require interpolation from all nodal values within that element.

An intermediate between these two cases is that of half-closed nodes, where $N_b$ nodes are placed only on some boundaries of each element. Specifically half-closed nodes satisfy the constraint that on each inter-element boundary between two neighbouring elements, at least one of the elements must have $N_b$ nodes placed on that boundary. The idea behind this is that it may allow for the retention of some of the advantages of using open or closed nodes, without incurring too much of the cost associated with strictly using one or the other. More details on how to determine which element places the $N_b$ nodes on each inter-element boundary are given in Sect. \ref{sect:halfclosed}.

\begin{figure}
    \centering
    \includegraphics[scale=0.4]{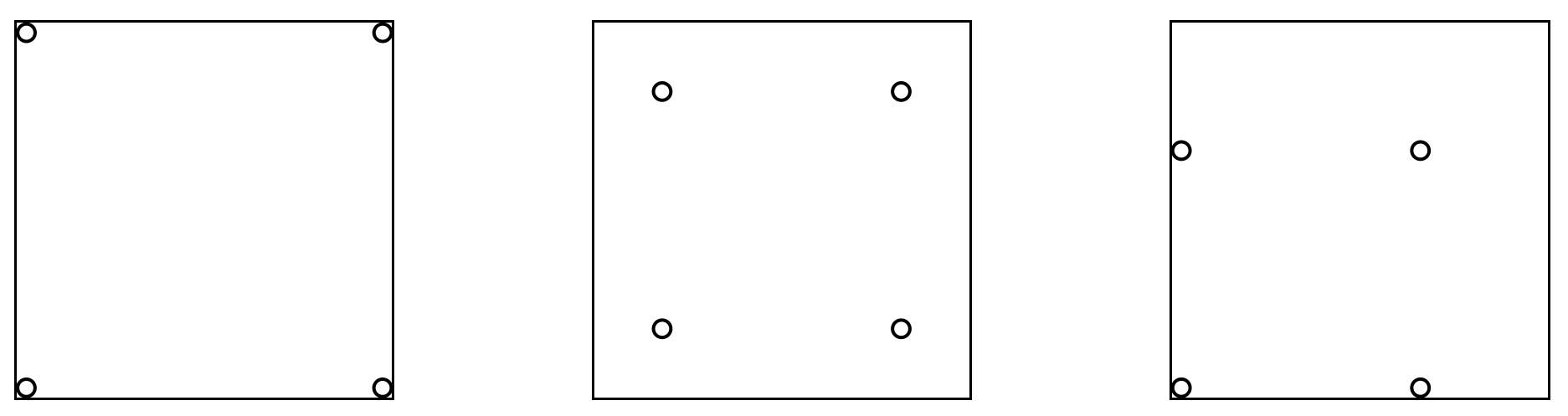}
    \caption{From left to right, examples of $p=2$ closed, open and half-closed nodes on a 2D quadrilateral element.}
    \label{fig:node_example}
\end{figure}

A schematic comparing closed, open, and half-closed nodes is shown in Fig. \ref{fig:node_example}. Some examples of the most popular nodes used in DG are Gauss-Lobatto, Chebyshev (both closed), Gauss-Legendre (open), and equidistant (either open or closed) nodes. Many of the above choices however are unique to quadrilateral elements however and are not easily extended to simplex elements.

\subsection{1st order equations}
We consider the first order scalar hyperbolic partial differential equation on the domain $\Omega$ with the inflow boundary or equivalently the Dirichlet boundary denoted as $\Gamma_{\text{in}}$
\begin{equation} \label{eq:conserv1d}
    \begin{aligned}
        \frac{\partial u}{\partial t} + \nabla &\cdot \boldsymbol{F}(u) = 0 \quad \text{in } \Omega \\
        u &= g \quad \text{on } \Gamma_{\text{in}}
    \end{aligned}
\end{equation}
where $\boldsymbol{F}(u)$ is the vector valued flux function. To find a solution $u \in V_h$, a weak formulation is obtained by multiplying Eq. (\ref{eq:conserv1d}) by a test function $v \in V_h$ and integrating by parts over each element $K_n \in \mathcal{T}_h$ to obtain
\begin{equation} \label{eq:weak1d_conserv}
    \frac{\partial}{\partial t}\int_{K_n} uv dx = \int_{K_n} \nabla v \cdot \boldsymbol{F}(u) dx - \int_{\partial K_n} v\hat{\boldsymbol{F}}(u^{\text{int}}, u^{\text{ext}}) \cdot \boldsymbol{n} ds \quad \forall{v} \in V_h
\end{equation}
where $u^{\text{int}}, u^{\text{ext}}$ denote the solution in the interior and exterior of the element $K_n$, and $\boldsymbol{n}$ is the outward pointing normal on the element boundary $\partial K_n$. The numerical flux $\hat{\boldsymbol{F}}$ is introduced only in the boundary term as boundary values of $\boldsymbol{F}$ are not required for computing the volume integral. For first order equations, many choices of numerical fluxes have been studied and used in various applications, and some common choices include the Lax-Friedrich, Roe, and Godunov fluxes.

To solve this system, the solution $u$ in $K_n$ is expressed as a linear combination of the basis functions $u = \sum_j \phi_j u_j, ~\phi_j \in V_h$, and the test function chosen to be an arbitrary basis function $\phi_i \in V_h$. Assuming that the numerical flux can be similarly represented as a linear combination of basis functions Eq. (\ref{eq:weak1d_conserv}) can be written in operator form as
\begin{equation} 
   M \frac{\partial}{\partial t} \boldsymbol{u} = \sum_d D^d \boldsymbol{F_d}
\end{equation}
where $\boldsymbol{u}, \boldsymbol{F_d}$ are vectors of size $|V_h|$, the number of basis functions, containing the values of the solution $u$ and the $d$-th component of the flux $F_d$ respectively. The matrix $M$ is known as the mass matrix and corresponds with the left side of Eq. (\ref{eq:weak1d_conserv}) while each matrix $D^d$ is $d$-th component of the discrete divergence operator acting on $F_d$ as
\begin{equation}
    \label{eq:div_full} D^dF_d = \int_{K_n} \frac{\partial v}{\partial x_d} F_d(u) ~dx - \int_{\partial K_n} v\hat{F}_d(u^{\text{int}}, u^{\text{ext}}) n_d ~ds \quad \forall{v} \in V_h
\end{equation}
Here $n_d, \hat{F}_d$ denote the $d$-th component of the outward facing normal on the boundary $\partial K_n$ and of the numerical flux $\hat{\boldsymbol{F}}$ respectively. The variable $x_d$ has also been introduced as the coordinate of the $d$-th dimension.

\subsection{2nd order operators} \label{sect:poi1d}
\begin{figure}
    \centering
    \includegraphics[scale=0.3]{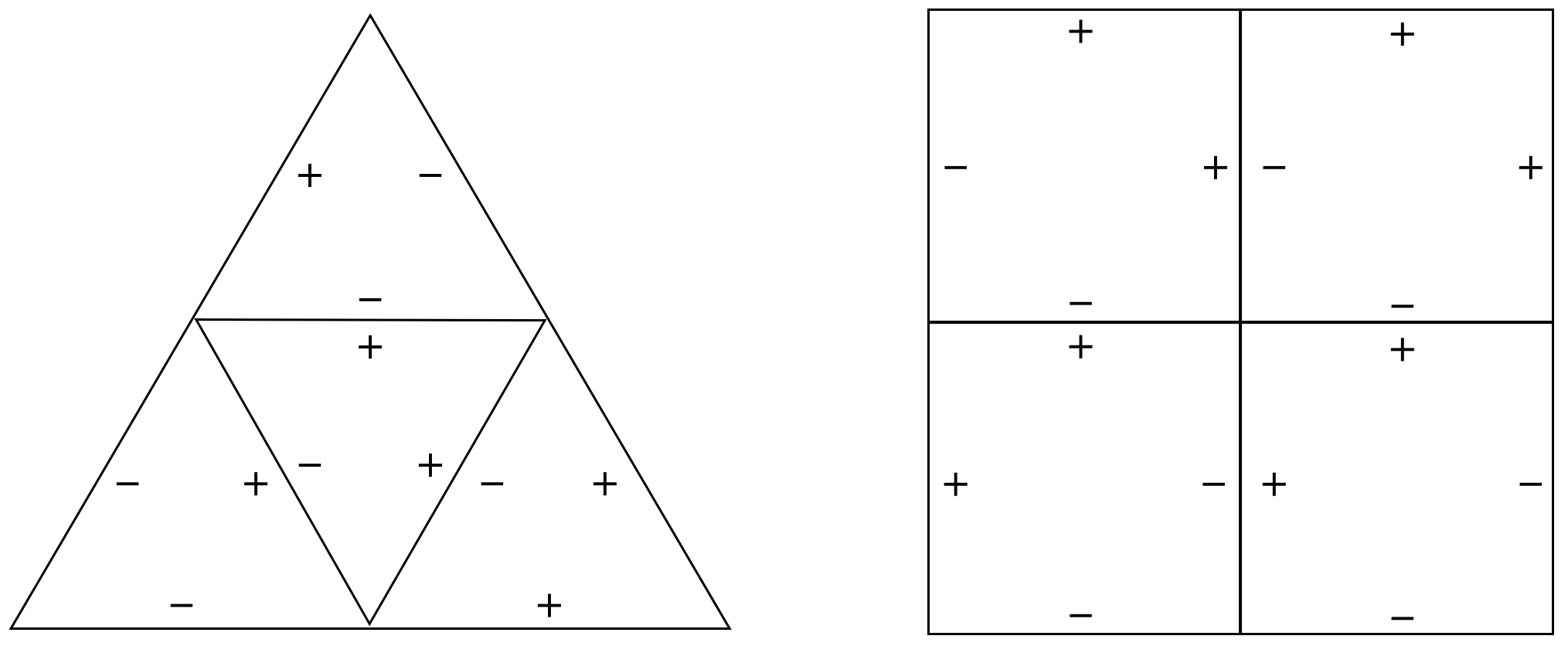}
    \caption{Example switch functions on a simplex and quadrilateral mesh in 2D. Switch function values are shown with $+,-$ for $+1,-1$ respectively. These switch functions are valid as the values along inter-element boundaries have opposing signs.}
    \label{fig:switch_example}
\end{figure}

For second order equations we consider Poisson's equation
\begin{equation} \label{eq:heat1d}
    \begin{aligned}
        - \Delta u &= f \quad \text{in } \Omega \\
        u &= g_D \quad \text{on } \Gamma_D \\
        \nabla u \cdot \bm{n} &= g_N \quad \text{on }\Gamma_N
    \end{aligned}
\end{equation}
as the model problem. To discretise this equation we use the LDG framework \cite{ldg1998}, which rewrites the equation into a split first order form. To do this a new variable is introduced for the gradient $\boldsymbol{q}=\nabla u$ and Eq. (\ref{eq:heat1d}) is rewritten to give the new system
\begin{align}
    - \nabla \cdot \boldsymbol{q} &= f \\
    \boldsymbol{q} &= \nabla u
\end{align}
A weak formulation is similarly obtained by multiplying this system of equations with test functions $\boldsymbol{\tau}, v$ and integrating by parts over each element $K_n$
\begin{align}
    \label{eq:ldg1} \int_{K_n} \bm{q} \cdot \nabla v ~ dx - \int_{\partial K_n} v \hat{\bm{q}} \cdot \bm{n} ~ds = \int_{K_n} fv ~ dx \quad &\forall v \in V \\
    \label{eq:ldg2} \int_{K_n} (\bm{q} - u \nabla) \cdot \bm{\tau} ~dx - \int_{\partial K_n} \hat{u} \bm{\tau} \cdot \bm{n} ~ds = 0 \quad &\forall \bm{\tau} \in V^d
\end{align}
In the LDG method fluxes $\hat{u}, \hat{\bm{q}}$ are specified using a so-called switch function defined on the element boundary $\partial K_n$. For the portion of $\partial K_n$ separating $K_n$ and some other element $K_m$, the switch function is denoted $S_n^m$ and given a value of either $-1$ or $+1$, where the only restriction is that $S_n^m + S_m^n = 0$. Examples of some valid switch functions are shown in Fig. \ref{fig:switch_example}. We consider strategies for assigning the switch function for both simplex and quadrilateral based elements in more detail in Sect. \ref{sect:halfclosed}.

Given a switch function the so-called minimal dissipation LDG flux $\hat{\bm{q}}$ can be defined as
\begin{equation} \label{eq:switch_q}
    \hat{\bm{q}} = 
    \begin{cases}
      \bm{q}|_{K_n} &\text{if } S_{n}^{m} < 0 \\
      \bm{q}|_{K_m}     &\text{if } S_{n}^{m} > 0
    \end{cases}
\end{equation}
and $\hat{u}$ similarly 
\begin{equation} \label{eq:switch_u}
    \hat{u} = 
    \begin{cases}
      u|_{K_n} &\text{if } S_{n}^{m} > 0 \\
      u|_{K_m}     &\text{if } S_{n}^{m} < 0
    \end{cases}
\end{equation}
This is sometimes known as an upwind-downwind flux, since the values of $\hat{u}, \hat{\boldsymbol{q}}$ are taken from opposing elements according to the switch.

On the domain boundary, the fluxes are chosen as
\begin{align}
    \hat{\bm{q}} &= 
    \begin{cases}
      \bm{q} - C_D(u - g_D)\bm{n}   &\text{on } \Gamma_D \\
      g_N \bm{n}              &\text{on } \Gamma_N \\
    \end{cases} \\    
    \hat{u} &= 
    \begin{cases}
      g_D &\text{on } \Gamma_D \\
      u   &\text{on } \Gamma_N
    \end{cases}
\end{align}
where $C_D \geq 0$ is a penalty parameter included for stabilisation on Dirichlet boundaries. Following the discussion in \cite{liu2021superconvergence}, this penalty is chosen to be non-zero only on domain boundaries where the switch function $S_n^m > 0$ is positive. It can be shown \cite{sherwin20062d} that the LDG method is stable with these flux definitions if the switch function satisfies the property that
\begin{equation} \label{eq:consistent_switch}
    \bigg| \sum_{K_m \in \mathcal{N}(K_n)} S_{n}^{m} ~\bigg| < |\mathcal{N}(K_n)|, ~\forall K_n \in \mathcal{T}_h
\end{equation}
where $\mathcal{N}(K_n)$ denotes the set of neighbouring elements of $K_n$. A switch function satisfying this property is termed a consistent switch function.

The system is solved similarly to the first order case by representing both $u, \boldsymbol{q}$ as linear combinations of basis functions in each element, and by picking the test functions to be each of the basis functions themselves. Through this Eqs. $(\ref{eq:ldg1})-(\ref{eq:ldg2})$ can be written in operator form as
\begin{equation}
    \begin{bmatrix}
        \boldsymbol{M} & G \\
        -D &
    \end{bmatrix}
    \begin{bmatrix}
        \boldsymbol{q} \\ u
    \end{bmatrix}
    =
    \begin{bmatrix}
        \boldsymbol{0} \\ Mf
    \end{bmatrix}
\end{equation}
where $\boldsymbol{M}$ is a stacked mass matrix, $D$ the discrete divergence, and $G$ the discrete gradient, where the gradient and divergence are related by the adjoint property $G = -D^T$. For example in two dimensions this can be expanded to give
\begin{equation}
    \label{eq:ldg_operator}
    \begin{bmatrix}
        M & & G^x \\
        & M & G^y \\
        -D^x & -D^y 
    \end{bmatrix}
    \begin{bmatrix}
        q^x \\ q^y \\ u
    \end{bmatrix}
    =
    \begin{bmatrix}
        0 \\ 0 \\ Mf
    \end{bmatrix}
\end{equation}
where $M$ is the standard mass matrix, and $D^d, G^d$ the components of the discrete divergence and gradient. In practice the solution $u$ can be solved for directly by taking the Schur complement of the above linear system without explicitly solving for $\bm{q}$ through the system
\begin{align}
    -Lu &= f \\
    \label{eq:ldg_laplacian} L = D^x M^{-1} &G^x + D^y M^{-1} G^y
\end{align}
where $L$ is known as the discrete LDG Laplacian.

\subsection{Spectral properties}
One key property of Galerkin methods is the independence of operator spectra with respect to basis representation for a given Finite Element function space $V_h$. In particular this property implies that eigenvalues of a linear operator are independent of the choice of basis representation of $V_h$ and therefore the choice of nodes. As such stability conditions dependent on the eigenvalues such as the CFL number remain unchanged irrespective of the type of nodes used. Furthermore it can be shown that symmetry of the discrete LDG Laplacian matrix is likewise maintained independent of the choice of nodes.

To show this we first define the change of basis matrix. Given two bases $\mathcal{B}, \tilde{\mathcal{B}}$ of the space $V_h$, the change of basis matrix is denoted $P_{\mathcal{B} \rightarrow \tilde{\mathcal{B}}}$ defined by 
\begin{equation}
    P_{\mathcal{B} \rightarrow \tilde{\mathcal{B}}} = \begin{bmatrix}
        [\phi_1]_{\tilde{\mathcal{B}}} & [\phi_2]_{\tilde{\mathcal{B}}} & ...
    \end{bmatrix}
\end{equation}
Here $[\phi_i]_{\tilde{\mathcal{B}}}$ denotes the coordinate of $\phi_i \in \mathcal{B}$ with respect to $\tilde{\mathcal{B}}$. For any linear operator $T:V_h \rightarrow V_h$ the matrix representations of $T$ with respect to these bases are related by a similarity transform
\begin{equation}
    \label{eq:change_basis}
    [T]_\mathcal{B} = P_{\mathcal{B} \rightarrow \tilde{\mathcal{B}}}^{-1} [T]_{\tilde{\mathcal{B}}} P_{\mathcal{B} \rightarrow \tilde{\mathcal{B}}}
\end{equation}
and therefore share the same eigenvalues. To apply this to the aforementioned DG operators we must first rearrange the equations. Taking the discrete gradient as an example we have
\begin{equation}
    Mq^d = -G^d u
\end{equation}
where $u, q^d$ are members of the same vector space. Thus the equation must be rearranged to 
\begin{equation}
    q^d = -\underbrace{M^{-1} G^d}_{\bar{G}^d} u
\end{equation}
so that $\bar{G}^d$, known as the nodal gradient operator, maps from $V_h$ to itself. Applying Eq. (\ref{eq:change_basis}) to this operator gives
\begin{equation}
    [\bar{G}^d]_\mathcal{B} = P_{\mathcal{B} \rightarrow \tilde{\mathcal{B}}}^{-1} [\bar{G}^d]_{\tilde{\mathcal{B}}} P_{\mathcal{B} \rightarrow \tilde{\mathcal{B}}}
\end{equation}
implying that the matrices for the nodal gradient under two different bases consequently have the same eigenspectra up to a change of basis. A similar argument can be constructed for the divergence and the same conclusion reached. 

For the LDG Laplacian we can plug definitions of the nodal gradient and nodal divergence into Eq. (\ref{eq:ldg_laplacian}) to get the nodal LDG Laplacian 
\begin{equation}
    \label{eq:nodal_laplacian}
    [M^{-1}L]_{\tilde{\mathcal{B}}} = [\bar{L}]_{\tilde{\mathcal{B}}} = \sum_d P_{\mathcal{B} \rightarrow \tilde{\mathcal{B}}} [\bar{D}^d\bar{G}^d]_{\mathcal{B}} P_{\mathcal{B} \rightarrow \tilde{\mathcal{B}}}^{-1}
\end{equation}
It also follows from the definition of the change of basis matrix that the mass matrices for two choices of bases are related by
\begin{equation}
    \label{eq:mass_basis_change}
    [M]_{\mathcal{B}} = P_{\mathcal{B} \rightarrow \tilde{\mathcal{B}}}^T [M]_{\tilde{\mathcal{B}}} P_{\mathcal{B} \rightarrow \tilde{\mathcal{B}}}
\end{equation}
Combining the two equations Eq. (\ref{eq:nodal_laplacian}) - (\ref{eq:mass_basis_change}) we get that
\begin{equation}
    [L]_{\tilde{\mathcal{B}}} = \sum_d P_{\mathcal{B} \rightarrow \tilde{\mathcal{B}}}^{-T} [D^dM^{-1}G^d]_{\mathcal{B}} P_{\mathcal{B} \rightarrow \tilde{\mathcal{B}}}^{-1} = P_{\mathcal{B} \rightarrow \tilde{\mathcal{B}}}^{-T} [L]_{\mathcal{B}} P_{\mathcal{B} \rightarrow \tilde{\mathcal{B}}}^{-1}
\end{equation}
establishing the claim that symmetry of the LDG Laplacian matrix is maintained under a change of basis.

\section{Sparsity patterns of the operators} \label{sect:sparsity}
In this section we investigate the sparsity of the matrix operators described in the previous section. To do so we use the notation $K_n \sim K_m$ to denote neighbouring elements $K_n, K_m$. We examine the sparsity pattern of three operators: the mass matrix $M$, the discrete divergence $D^d$, and the discrete LDG Laplacian $L$. In particular we are interested in how the sparsity pattern of each operator changes when using each different type of node (open, closed, half-closed) described in Sect. \ref{sect:nodes}.

\subsection{Mass matrix} \label{sect:mass}
The mass matrix has entries defined by
\begin{equation}
    M_{ij} = \sum_{n} \int_{K_n} \phi_i \phi_j dx
\end{equation}
which contains only volume integral terms, meaning it is a block diagonal operator for any choice of nodes. This is because the product of any two basis functions $\phi_i \phi_j \neq 0$ only if they belong to the same element. Each block diagonal is in general a dense block of size $\mathcal{O}(p^d)$, where $p$ is the polynomial degree and $d$ the spatial dimension.

For some special choice of nodes the mass matrix can be reduced simply to a diagonal operator if
\begin{equation} \label{eq:mass_delta}
    \sum_{n} \int_{K_n} \phi_i \phi_j dx = C_i \delta_{ij}
\end{equation}
where $C_i = \sum_{n} \int_{K_n} \phi_i^2 dx$, that is if the basis functions are mutually orthogonal. For quadrilaterals, for any polynomial degree $p$ examples of open nodes with this property include the Gauss-Legendre nodes and of half-closed nodes the Gauss-Radau nodes. There however do not exist any closed nodes on quadrilaterals with this property, although it is common in for instance the DG-SEM method to accept a slight under-integration and use approximate diagonal mass matrices with Gauss-Lobatto nodes. For simplices, to the authors' knowledge, there are as of yet no known nodes of any type that produce a diagonal mass matrix for arbitrary polynomial degrees $p$.

\subsection{Divergence/gradient operators}
We examine the sparsity pattern for the divergence defined in Eq. (\ref{eq:div_full}). While the exact form of the divergence depends on the numerical flux used, to study the sparsity pattern we can without loss of generality take the entries of the discrete divergence to be
\begin{equation} \label{eq:div_nd}
    D^d_{ij} = \sum_n \int_{K_n} \frac{\partial \phi_i}{\partial x_d} \phi_j dx - \int_{\partial K_n} \phi_i \phi_j n_d ds
\end{equation}
The sparsity pattern of the gradient operator is also determined here given they are related via a simple transpose. As with the mass matrix, the divergence operator has for each element $K_n$ a dense block on the diagonal corresponding to the volume integral term. In addition, due to the boundary integral term the divergence operator will in general also have non-zero entries in off-diagonal blocks of the matrix coupling neighbouring elements on an element boundary $\partial K_n$.

\begin{figure}
    \centering
    \includegraphics[scale=0.33]{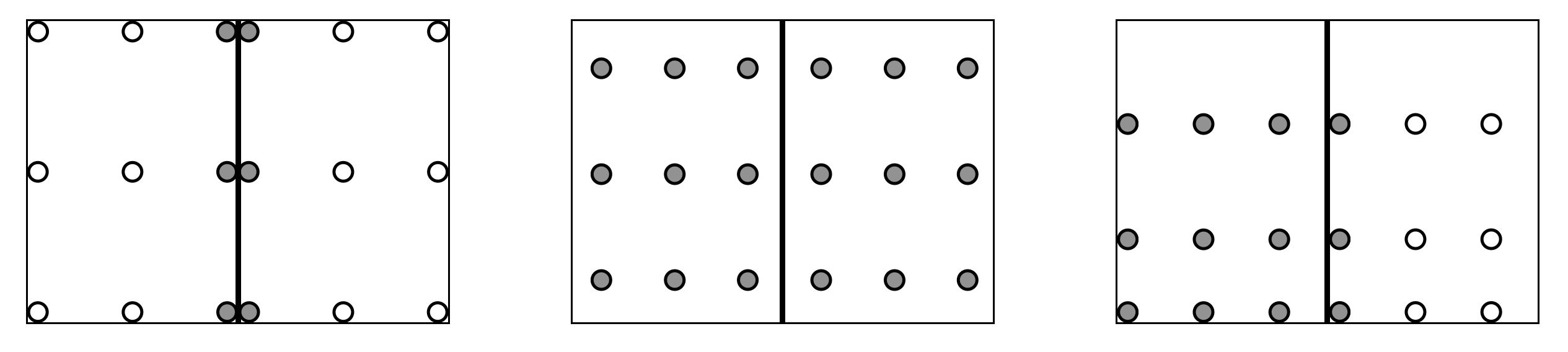}
    \caption{Communication pattern for boundary terms between neighbouring elements in discrete divergence operator for different nodes. On the left for closed nodes only the nodal values on the boundary are needed to evaluate the boundary terms. In the middle for open nodes as no nodes are on the boundary all the nodal values in both elements are needed from both elements. On the right for half-closed nodes, in the left element only the nodal values on the boundary are needed whilst on the right as no nodes are on the boundary all the nodal values in the element are needed similar to the open case.}
    \label{fig:sparsity_div}
\end{figure}

\begin{figure}
    \centering
    \includegraphics[scale=0.33]{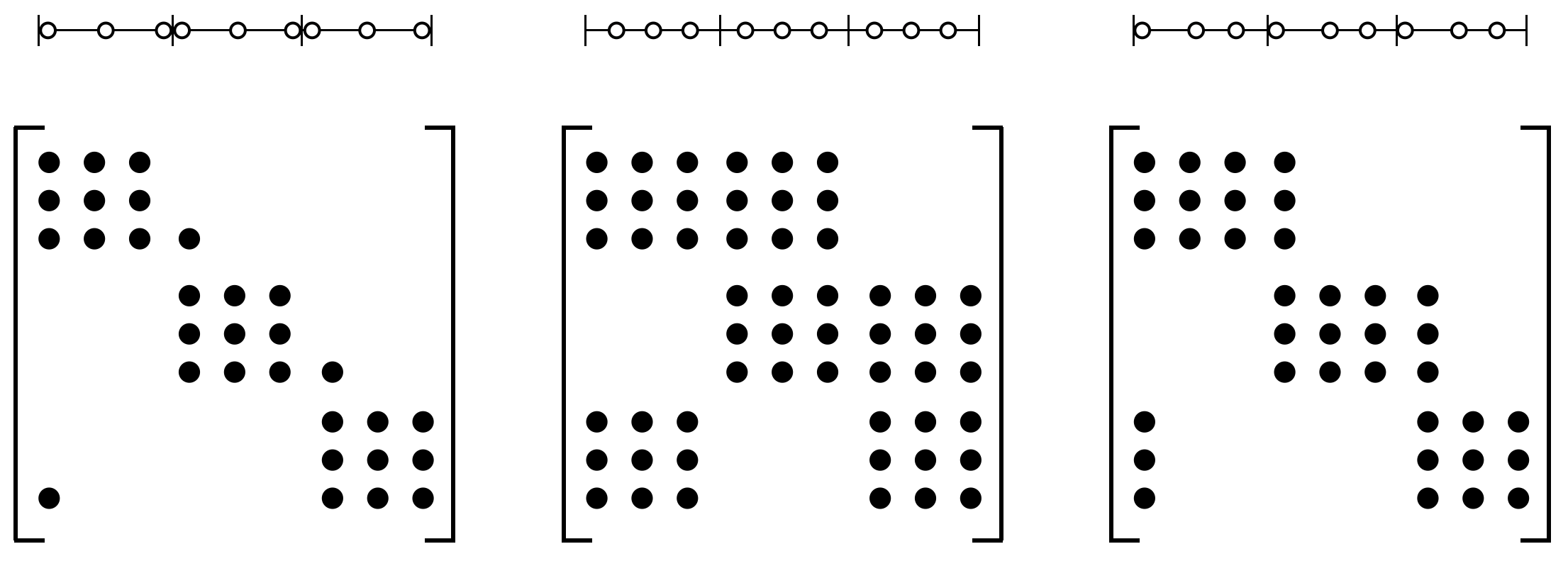}
    \caption{Comparison sparsity of the divergence operator on 1D mesh with three $p=2$ elements and periodic boundary conditions. The flux here at every inter-element boundary is taken to be $\hat{F}(F_L, F_R) = F_R$, the value from the right element. From left to right, the sparsity pattern of the divergence with this numerical flux using closed, open, and half-closed nodes is shown.}
    \label{fig:sparsity_div_comp}
\end{figure}

The sparsity patterns of the off-diagonal blocks in the divergence operator for each of the three types of nodes are considered separately. The communication patterns for all three types of nodes for the boundary terms of the discrete divergence operator are shown in Fig. \ref{fig:sparsity_div}. A side-by-side comparison of the sparsity pattern of the divergence for the three types of nodes on an example 1D mesh is also shown in Fig. \ref{fig:sparsity_div_comp}.

\subsubsection{Open nodes}
The simplest case to consider is that of open nodes as no distinction needs to be made between nodes interior to an element and those on the boundary. In this case evaluation of function values on the boundary between $K_n \sim K_m$ requires interpolation of all nodal values from elements $K_n$ and $K_m$. This in general results in a communication pattern where all nodes of elements $K_n, K_m$ communicate with one another for the boundary integral term and implies that off-diagonal blocks in the divergence/gradient for neighbouring elements $K_n \sim K_m$ will be fully dense, with each having $\mathcal{O}(p^{2d})$ entries.

\subsubsection{Closed nodes}
For closed nodes, $N_b$ nodes are placed by each element on all its boundaries such that all basis functions corresponding to interior nodes are zero on the boundaries. As a result the boundary terms in Eqs. (\ref{eq:ldg1}) - (\ref{eq:ldg2}) can be ignored for the interior nodes for both the divergence and gradient operators. This implies that in off-diagonal blocks corresponding to communication between neighbouring elements $K_n \sim K_m$, rows and columns corresponding to interior nodes of either $K_n$ or $K_m$ will not contain any non-zero entries.

Only the $2N_b$ boundary nodes from each element on the boundary between $K_n \sim K_m$ will communicate with one another for the boundary integral term in the discrete divergence and gradient. This is due to the fact that with closed nodes the value of a function on an element boundary can always be determined using only the $2N_b$ nodal values at the boundary. Therefore the off-diagonal block in the divergence/gradient operator for neighbouring elements $K_n \sim K_m$ will each contain $\mathcal{O}(N_b^2)$ entries, or equivalently $\mathcal{O}(p^{2(d-1)})$ entries.

\subsubsection{Half-closed nodes}
In between open and closed is the case of half-closed nodes where on each inter-element boundary between two elements, $N_b$ nodes are required to be placed on that boundary by at least one of the elements. This implies that on any inter-element boundary, the basis function values on that boundary for at least one element can be inferred solely using nodal values at the boundary, whilst for the other interpolation from all its nodal values may be required as with the open case. Thus the off diagonal block for neighbouring elements $K_n \sim K_m$ has at most $\mathcal{O}(p^{2d-1})$ entries, in between that of the open and closed cases.

\subsection{LDG Laplacian} \label{sect:sparsity_ldg}
The LDG Laplacian in Eq. (\ref{eq:ldg_laplacian}) can be expressed in terms of the mass and divergence operators as
\begin{equation} \label{eq:LDG_laplacian}
    L = \sum_d D^d M^{-1} G^d
\end{equation}
where now the fluxes in the divergence $D^d$ and similarly in the gradient $G^d$ are defined using a switch function. The only point of difference here in the gradient/divergence operators from the previous subsection is in the flux which for the divergence is now specified in Eq. (\ref{eq:switch_q}) and for the gradient in Eq. (\ref{eq:switch_u}). This implies that for neighbouring elements $K_n \sim K_m$, nodes from an element $K_n$ will communicate only with nodes from $K_m$ in the divergence operator if the switch $S_n^m > 0$. For the gradient the reverse is true due to the upwind-downwind flux structure of LDG, so that nodes from element $K_n$ will communicate only with those from $K_m$ if the switch $S_n^m < 0$. 

\begin{figure}
    \centering
    \includegraphics[scale=0.4]{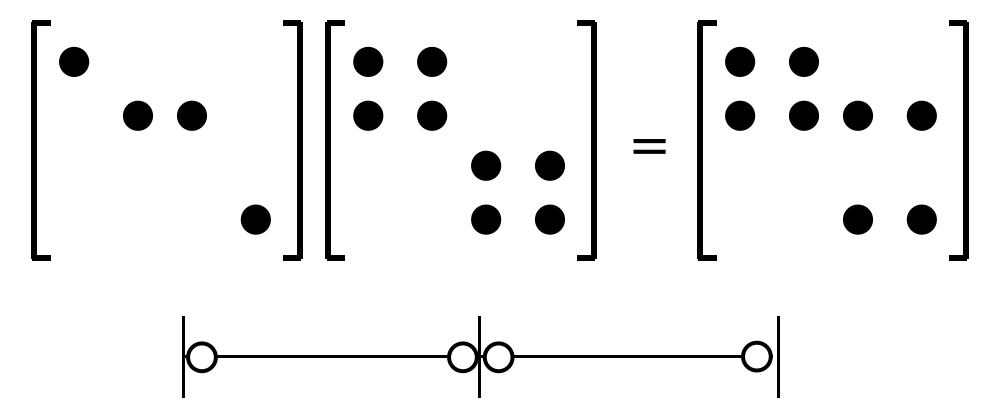}
    \caption{Effect of applying an inverse mass matrix (center matrix) on the right of an arbitrary operator (left matrix). In this example the operators correspond to the two element $p=1$ shown at the bottom. The operator initially has only one non-zero per column. After applying the mass matrix on the right this spreads the non-zero entries across the rows of each block (right matrix).}
    \label{fig:sparsity_massprod}
\end{figure}

From the above discussion, the mass matrix is block diagonal where each block corresponds to the nodes of a single element. Furthermore the blocks are in general dense except for some special choices of nodes described in Sect. \ref{sect:mass}. As such when applied to another operator on the right, the inverse mass matrix in general spreads any non-zero entry in the $(i,j)$ position of the operator across a row to all columns corresponding to nodes in the same element as the $i$-th node. This effect is illustrated in Fig. \ref{fig:sparsity_massprod}.

Since the LDG Laplacian is the sum of products of these operators, we are able to examine the sparsity patterns of the Laplacian for each type of the three nodes. We focus on the sparsity of the off-diagonal blocks as the diagonal blocks are fully dense in each of the three cases. As before each type of node is considered separately. A side-by-side comparison of the sparsity pattern of the LDG Laplacian for the three types on an example 1D mesh is also shown in Fig. \ref{fig:sparsity_ldg_comp}.

\subsubsection{Open nodes}
According to the switch function, for the divergence all nodes in an element $K_n$ will communicate with those from a neighbouring element $K_n \sim K_m$ if $S_n^m > 0$. For the gradient operator the opposite is true, in that all nodes from neighbouring elements $K_n \sim K_m$ will communicate if $S_n^m < 0$. As discussed in the previous subsection for open nodes, these off-diagonal blocks are fully dense.

This means that right multiplication by the inverse mass matrix does not affect the sparsity pattern of the divergence operator in this case. The pattern can therefore be deduced by only considering the pattern of the product of the divergence with the gradient. This implies that the Laplacian for open nodes will have non-zero entries in off-diagonal blocks where each node in an element $K_n$ denoted $s_j^n$ communicates with:
\begin{enumerate}
    \item All nodes $\{ s_j^m \}$ from neighbour elements $K_n \sim K_m$
    \item All nodes $\{ s_j^o \}$ from second neighbour elements $K_n \sim K_m \sim K_o$ if:
    \begin{enumerate}
        \item $S_n^m > 0$ and,
        \item $S_m^o < 0$
    \end{enumerate}
\end{enumerate}
\begin{figure}
    \centering
    \includegraphics[scale=0.35]{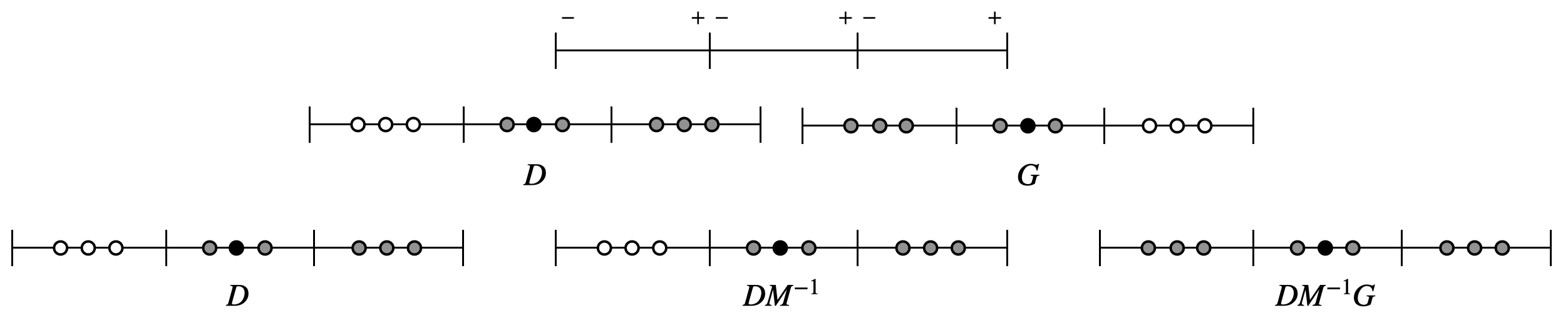}
    \caption{Communication pattern for open nodes in LDG operators on a 1D $p=2$ mesh with switch function shown at the top. Nodes in black correspond to a node $i$, and shaded nodes all nodes $j$ where the $(i,j)$ entry of the relevant matrix is non-zero. Middle row shows divergence and gradient communication patterns for node shaded in black. On the bottom row, the operators are multiplied to give the communication pattern of the Laplacian for the node in black.}
    \label{fig:sparsity_ldg_o}
\end{figure}
The communication on an 1D example is shown in Fig. \ref{fig:sparsity_ldg_o}. Fig. \ref{fig:sparsity_ldg_comp} also shows the resulting communication and sparsity patterns on a simple 2D example.

\subsubsection{Closed nodes}
The only difference in the sparsity patterns of the divergence and gradient operators with closed nodes compared with open ones is in the sparsity of the off-diagonal blocks. In the case where closed nodes are used, as described in the previous subsection the off-diagonal blocks are sparse and only contain non-zero entries corresponding to communication between the $2N_b$ nodes on the boundary between neighbouring elements $K_n \sim K_m$. This is in contrast to the off-diagonal blocks when using open nodes which are fully dense. 

\begin{figure}
    \centering
    \includegraphics[scale=0.35]{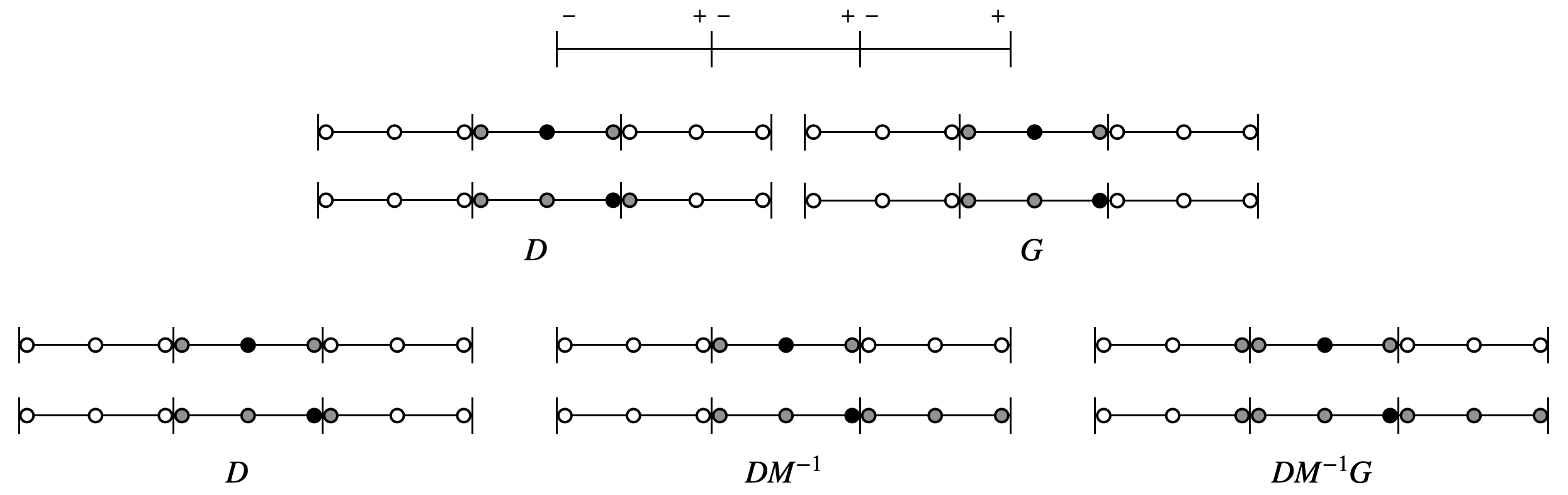}
    \caption{Communication pattern for closed nodes in LDG operators on a 1D $p=2$ mesh for switch function shown at the top. Nodes in black correspond to a node $i$, and shaded nodes all nodes $j$ where the $(i,j)$ entry of the relevant matrix is non-zero. Middle rows shows divergence and gradient communication patterns for interior and boundary nodes shaded in black. On the bottom row, the operators are multiplied to give the communication pattern of the Laplacian for the nodes in black.}
    \label{fig:sparsity_ldg_c}
\end{figure}

However to compute the Laplacian matrix the divergence operator is first multiplied on the right by the inverse mass matrix. This multiplication on the right of the mass inverse to the divergence has the effect of spreading non-zero entries as shown in Fig. \ref{fig:sparsity_massprod}. This is then multiplied on the right by the gradient operator to form the LDG Laplacian. Overall this gives that the Laplacian with closed nodes has non-zero entries corresponding to communication by nodes $s_j^n$ in element $K_n$ with:
\begin{enumerate}
    \item All boundary nodes $\{ s_j^m \}$ on the boundary between neighbour elements $K_n \sim K_m$ where $S_n^m < 0$
    \item All nodes $\{ s_j^m \}$ from neighbour element $K_n \sim K_m$ if:
    \begin{enumerate}
        \item $S_n^m > 0$ and,
        \item the node $s_j^n$ is on the boundary between $K_n, K_m$
    \end{enumerate}
    \item All boundary nodes $\{ s_j^o \}$ on boundary between $K_m \sim K_o$ where $K_o$ is a second neighbour element $K_n \sim K_m \sim K_o$ if:
    \begin{enumerate}
        \item $S_n^m > 0$ and,
        \item $S_m^o < 0$ and,
        \item the node $s_j^n$ is on the boundary between $K_n, K_m$
    \end{enumerate} 
\end{enumerate}
The communication on an 1D example is shown in Fig. \ref{fig:sparsity_ldg_c}. Fig. \ref{fig:sparsity_ldg_comp} also shows the communication pattern and resulting sparsity on a simple 2D example.

\subsubsection{Half-closed nodes}
We make the observation that the differences in sparsity pattern of the Laplacian in the open and closed cases can be attributed simply to having $N_b$ boundary nodes on the boundary $K_n \sim K_m$ in element $K_n$ where the switch function satisfies $S_n^m > 0$. This suggests the following strategy for placing half-closed nodes: on an inter-element boundary between two elements $K_n \sim K_m$, element $K_n$ must place $N_b$ nodes on that boundary if $S_n^m > 0$, and vice versa. As the switch function is constrained to satisfy $S_n^m + S_m^n = 0$, this guarantees at least one element on each inter-element boundary will place $N_b$ nodes on that boundary. If half-closed nodes are placed in this fashion, despite that the divergence and gradient operators have a slightly increased number of entries in off-diagonal blocks than in the closed case, the resulting sparsity pattern of the LDG Laplacian for half-closed nodes will be identical to that of closed nodes. 

\begin{figure}
    \centering
    \includegraphics[scale=0.35]{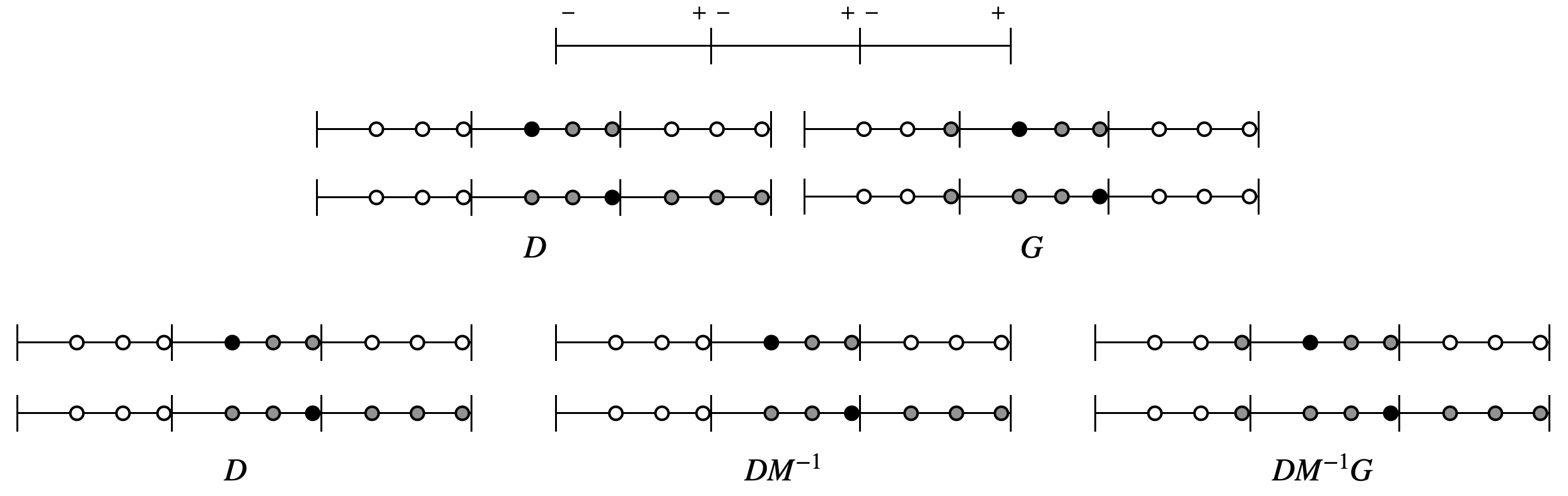}
    \caption{Communication pattern for half-closed nodes in LDG operators on a 1D $p=2$ mesh for switch function shown at the top. Nodes in black correspond to a node $i$, and shaded nodes all nodes $j$ where the $(i,j)$ entry of the relevant matrix is non-zero. Middle rows shows divergence and gradient communication patterns for interior and boundary nodes shaded in black. On the bottom row, the operators are multiplied to give the communication pattern of the Laplacian for the nodes in black.}
    \label{fig:sparsity_ldg_h}
\end{figure}

\begin{figure}
    \centering
    \includegraphics[scale=0.35]{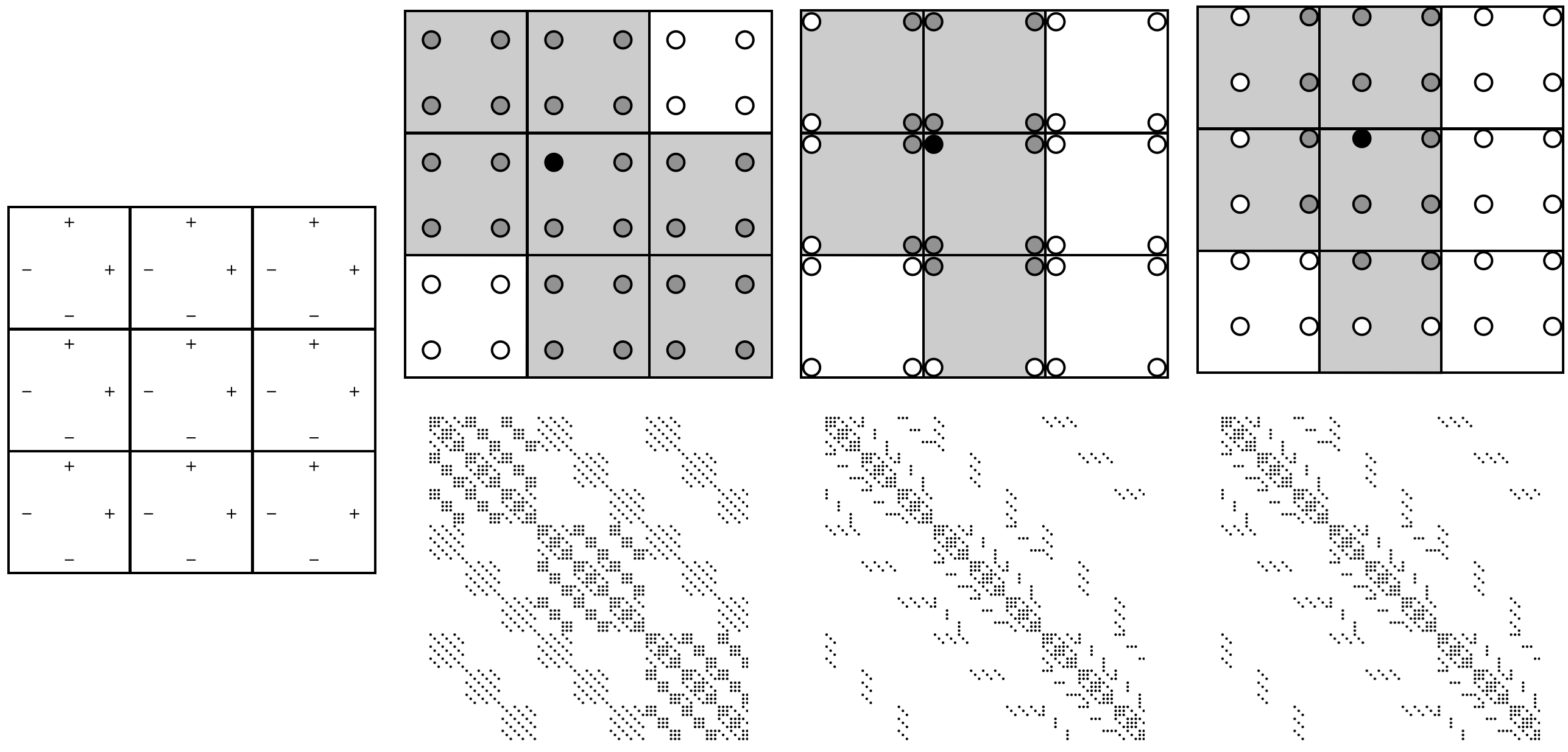}
    \caption{Comparison of LDG communication pattern and sparsity of the LDG Laplacian on a 3x3 Cartesian mesh using $p=1$ elements in 2D. Nodes shaded in black correspond to a node $i$, and shaded nodes all nodes $j$ where the $(i, j)$ entry of the relevant matrix is non-zero. Elements with at least one node included in the communication pattern are shaded. From left to right the sparsity pattern using open, closed, and half-closed nodes are shown. The communication pattern and therefore sparsity for the closed and half-closed cases are identical, with each having 729 non-zero entries. The open case suffers from a decrease in sparsity, with 1377 non-zero entries.}
    \label{fig:sparsity_ldg_comp}
\end{figure}

This conclusion holds for all element shapes in all dimensions. This is because no assumptions have been made up to this point on the type of element considered, and the only requirement has been the existence of a valid switch function on a given mesh. An example of the communication pattern is shown for half-closed nodes in 1D in Fig. \ref{fig:sparsity_ldg_h} and the resultant communiction an sparsity patterns on a 2D example in Fig. \ref{fig:sparsity_ldg_comp}.

\section{Assembly cost}
\subsection{Volume terms}
Once the mass matrix and divergence operators have been obtained, this is sufficient to construct all the operators outlined in the previous section. This is as the gradient operator is related to the divergence operator via a simple transpose, and the LDG Laplacian in Eq. (\ref{eq:LDG_laplacian}) is constructed via a composition of the gradient, divergence and mass matrices.

To assemble the mass and divergence matrices, a volume integral on the element $K_n$ of the form
\begin{equation}
    P_i = \int_{K_n} \phi_i p(x) ~dx
\end{equation}
needs to be computed, where $\phi_i$ denotes the $i$-th basis function, and $p(x)$ some polynomial of degree of degree less or equal to that of $\phi_i$. For the mass matrix, the polynomial under consideration is $p(x) = \phi_j(x)$, whilst for the divergence matrix $p(x) = \frac{\partial}{\partial x^d} \phi_j$. As $p(x)$ is of degree less than equal to that of $\phi_i$, we can expand it as
\begin{equation}
    p(x) = \sum_j p_j \phi_j(x)
\end{equation}
and so the integral can be computed as
\begin{equation}
    \begin{split}
    P_i &= \sum_j \int_{K_n} \phi_i \phi_j ~dx \cdot p_j  \\
    &= \sum_j M_{ij} p_j 
    \end{split}
\end{equation}
implying that the cost of assembling the volume term is equal to the cost of left matrix multiplication with the mass matrix. It is therefore desirable for the mass matrix to be as sparse as possible to minimise the assembly cost.

As discussed in the previous section, DG mass matrices are block diagonal where the number of blocks is equal to the number of elements $|\mathcal{T}_h|$. These blocks are in general fully dense, except when the basis functions are pairwise orthogonal satisfying Eq. (\ref{eq:mass_delta}) in which case the mass matrix is strictly diagonal. This implies that the efficiency of a nodal DG method in terms of both sparsity and assembly cost is optimised if half-closed or closed nodes that produce pairwise orthogonal basis functions are chosen for the discretisation.

\subsection{Boundary terms}
For the boundary integral terms in Eq. (\ref{eq:div_nd}) we are interested in computing integrals of the form
\begin{equation}
    Q_i = \int_{\partial K_n} \phi_i q(x) ~ds
\end{equation}
where $\partial K_n$ denotes the boundary of an element $K_n$, and $q(x)$ some polynomial of degree less than or equal to that of $\phi_i$. A similar procedure can then be performed as above to give
\begin{equation}
    \begin{split}
    Q_i &= \sum_j \int_{\partial K_n} \phi_i \phi_j ~ds \cdot q_j  \\
    &= \sum_j S_{ij} q_j 
    \end{split}
\end{equation}
where $S_{ij}$ is a DG `boundary' mass matrix of sorts. The cost of assembling the boundary terms thus depends on the sparsity of this operator $S$, which is also block diagonal but in general also much sparser than the standard mass matrix. Its sparsity can similarly be reduced to a strictly diagonal matrix in the case that
\begin{equation}
    \int_{\partial K_n} \phi_i \phi_j ~ds = \delta_{ij}
\end{equation}
that is if basis functions are pairwise orthogonal with respect to boundary integration.

\section{Half-closed discontinuous Galerkin methods} \label{sect:halfclosed}
As shown in Sect. \ref{sect:sparsity}, half-closed and closed nodes produce operators with increased sparsity compared to those using open nodes. For first order operators, discretisations with half-closed nodes have slightly denser off-diagonal blocks than when using closed nodes. However if half-closed nodes are placed carefully the sparsity of second order LDG Laplace operators is identical for discretisations using half-closed or closed nodes. Specifically this is achieved if for a given mesh with an assigned switch function, $N_b$ nodes are placed on boundaries of an element $K_n$ where the switch function $S_n^m > 0$.

The advantage of using half-closed nodes over closed nodes comes from the extra freedom with which nodes can be placed on each element, without the constraint of needing exactly $N_b$ nodes on every element boundary unlike with closed nodes. For quadrilateral elements, this extra degree of freedom allows nodes to be placed at Gauss-Radau quadrature points, which produces pairwise orthogonal basis functions. This implies that the mass matrix using these nodes to be strictly diagonal, which following the discussion in the previous section also allows for more efficient operator assembly.

\subsection{Quadrilateral elements}

\subsubsection{Gauss-Radau nodes}
In 1-dimension with $n$ points on some interval, it is well known that the highest possible quadrature precision is attained through using the Gauss-Legendre points which is able to integrate exactly polynomials of degree up to $2n-1$. The Gauss-Legendre points are however open since none of the points lie on the boundary. If the restriction is added so that nodes lie on both boundaries of the interval, that is if they are to be closed, the highest possible quadrature precision is attained through using the Gauss-Lobatto points, which integrates polynomials up to degree $2n-3$ exactly.

If a looser restriction of only being half-closed is instead allowed, the highest possible quadrature precision is $2p-2$ which is achieved by the Gauss-Radau points. This implies that for a given polynomial degree $p$, if the $p+1$ half-closed nodes in each element in 1D are placed at the Gauss-Radau points, the integral terms in each of the operators in Sect. \ref{sect:sparsity} can be computed exactly using only the nodal values of each element without requiring the usual assembly procedure of interpolating nodal values to ones at separate quadrature points. This in particular implies that in this case the mass matrix when using Gauss-Radau nodes is in fact diagonal, as opposed to only being block diagonal.

\begin{figure}
    \centering
    \includegraphics[scale=0.35]{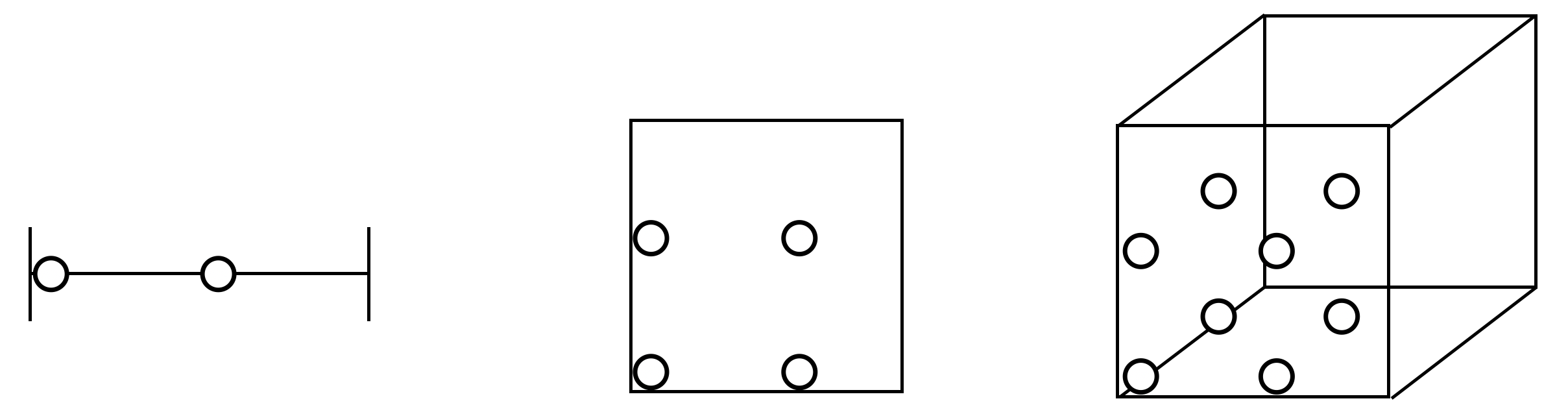}
    \caption{$p=1$ half-closed nodes on quadrilateral elements in 1D, 2D and 3D. Nodes are placed on exactly half of the boundaries of each element.}
    \label{fig:gaussradau}
\end{figure}

In higher dimensions, the half-closed nodes can be simply taken to be the $d$-dimensional outer product of the Gauss-Radau nodes in 1-dimension. This is shown in Fig. \ref{fig:gaussradau}. Similar to 1-dimension, these nodes can be used as quadrature points to integrate $d$-dimensional outer products of degree at most $2p-2$ polynomials exactly. This in particular implies that the mass matrix will remain diagonal as with the 1-dimensional case. Moreover if the boundaries of an element $K_n$ are linear, the nodal values can be used also to calculate both the volume integral terms and the boundary integral terms without requiring any interpolation to separate boundary quadrature points. 

One thing to note when using Gauss-Radau nodes is that for opposing boundaries in an element only one of them will have $N_b$ nodes placed on it. This is a consequence of the 1-dimensional Gauss-Radau points having a boundary point on only one of the endpoints. Given that in our methodology boundary nodes in the half-closed nodes are placed according to a given switch function we therefore require a way to assign switch functions on quadrilateral meshes that satisfy the property of having opposite values for opposing boundaries.

\subsubsection{Switch function}
To assign a valid switch function on a quadrilateral mesh, we adopt the simple algorithm used in \cite{linedg}. This algorithm is constructed using the fact that switch functions on quadrilateral elements should satisfy the two properties:
\begin{enumerate}
    \item The switch on a shared boundary between neighbouring elements $K_n \sim K_m$ should have opposing signs, that is $S_n^m + S_m^n = 0$
    \item The switch on opposing boundaries in the same element should have opposing signs
\end{enumerate}
These two properties are always possible to satisfy on an quadrilateral mesh in the absence of hanging nodes. This is because if an element boundary is selected at random and a global line drawn across the mesh connecting opposing element boundaries, this line is guaranteed to either loop around on itself or connect two domain boundaries. Switch functions on quadrilateral elements can therefore be generated using a simple iterative algorithm, where at each step an arbitrary boundary between two neighbouring elements $S_n \sim S_m$ is picked and the switch arbitrarily assigned as $S_n^m = -1$. The two properties listed above then propagate the switch in an alternating pattern along a sequence of boundaries in both directions until either a domain boundary is reached or the sequence repeats. This procedure can then be iterated until all boundaries have been processed. 

\begin{figure}
    \centering
    \includegraphics[scale=0.35]{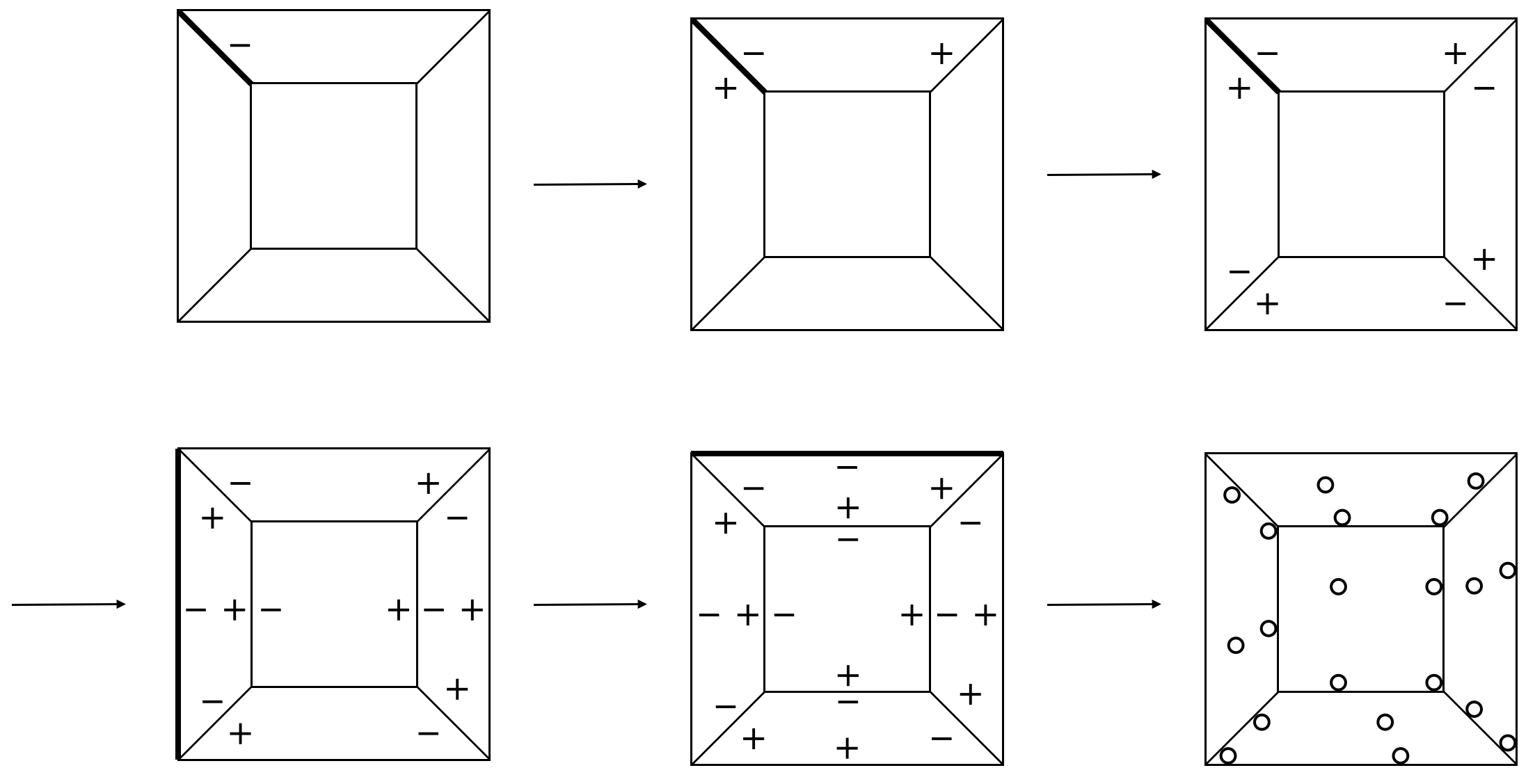}
    \caption{Schematic of switch function and node assignment on quadrilateral meshes. At each step an interior boundary between two elements shown bolded is picked and a switch of $-1$ arbitrarily assigned to one of the elements. This switch is then propagated in alternating fashion and the process repeated until all switches have been assigned. Finally half-closed nodes are placed in each element according to the switch function with $N_b$ boundary nodes placed where the switch has value $+1$.}
    \label{fig:quad_switch}
\end{figure}

Clearly, switch functions generated via this procedure are guaranteed to satisfy the two above properties and a schematic of this procedure on a mesh in 2-dimensions is shown in Fig. \ref{fig:quad_switch}. This is a simple consequence of the fact that at each step of the algorithm the switch is propagated by using the two properties themselves. Once a switch function has been chosen on a mesh in this way, Gauss-Radau nodes can be placed on each element $K_n$ with boundary nodes placed according to where the switch function $S_n^m > 0$. 

\subsubsection{Boundary quadrature}
The Gauss-Radau points can also be used as quadrature points for computing the boundary integral terms in the divergence operator. This is as the Gauss-Radau points achieve the same quadrature precision on element boundaries, which are similarly quadrilaterals but of one lower dimension. We outline briefly an efficient way to employ these as quadrature points, such that the cost of using them is the same as when using closed nodes.

On every inter-element boundary between two quadrilateral elements $K_n \sim K_m$, it is guaranteed that only one of the elements has $N_b$ nodes placed on the boundary. Assuming $K_n$ to be the element with nodes placed on the boundary, evaluations of the solution in the element $K_n$ on the boundary can be computed cheaply via only interpolation using its boundary nodes. However functions evaluations on this boundary in $K_m$ are expensive unless the evaluation points considered happen to align up with the Gauss-Radau points in the element so that they can be calculated via interpolation along a line. 

\begin{figure}[h]
    \centering
    \includegraphics[scale=0.4]{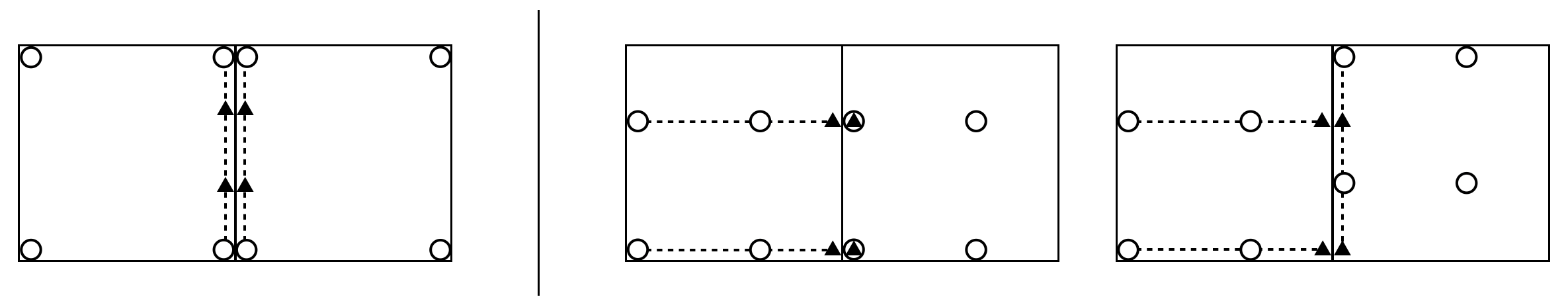}
    \caption{Schematic of boundary quadrature rules for half-closed and closed nodes in 2D. Nodes are shown as dots, while quadrature points are shown using triangles. On the left, boundary quadrature with closed nodes in shown, where 1D line interpolation between the boundary nodes is used to evaluate the solution at the quadrature points. On the right, if the boundary quadrature points are picked to always align with those in the element without nodes on the boundary, 1D interpolation can similarly be used to compute the integral.}
    \label{fig:boundary_quadrature}
\end{figure}

This implies that for efficient boundary quadrature computation, the points chosen should be the boundary Gauss-Radau points aligned with those from $K_m$. This is illustrated in Fig. \ref{fig:boundary_quadrature}. Performed in this way, only line interpolation is required for boundary integration in $K_m$, and line (2D) or face (3D) interpolation needed for boundary integration in $K_n$. This is similar to the cost incurred computing the boundary terms when using closed nodes, where similar interpolation is performed.

\subsection{Simplex elements}
As described in Sect. \ref{sect:discretisation} the number of nodes required on a simplex to define a polynomial of degree $p$ is equal to $N={p+d \choose d}$. However for $p \geq 3$ in 2-dimensions, and for $p \geq 2$ in 3-dimensions the minimum number of quadrature points required to exactly integrate a degree $2p$ polynomial is greater than the number of nodes $N$ in a degree $p$ simplex element \cite{dunavant1985quadrature,witherden2015quadrature}. Furthermore quadrature points in these optimal quadrature rules do not have any points placed on element boundaries and as a result strictly give open nodes.

This implies that there are no analogues of the half-closed Gauss-Radau or the open Gauss-Legendre nodes for simplex elements. Unlike in quadrilaterals there unfortunately do not exist on simplices any nodal discretisation that can integrate the mass matrix exactly for arbitrarily high polynomial degrees $p$, unless extra interpolation is employed to a separate set of quadrature points. This means that the mass matrices produced on simplex elements will in general not be diagonal if polynomial basis functions are used. We note however that there are however some non-polynomial choices of basis functions on simplices that have been found that do produce diagonal mass matrices \cite{hicken2016multidimensional}, although these basis functions do not currently have a simple closed-form expression.

\section{Linear solvers}
In this section we discuss two types of linear solver technique that aim to take advantage of the sparsity structure of the DG divergence and Laplacian operators when using half-closed nodes. We note that the constructions here can be equally applied to closed nodes, but not to open ones. This follows from the discussion of the sparsity of the operators in Sect. \ref{sect:sparsity}, wherein similar sparsity patterns for these operators are obtained when using closed or half-closed nodes. Specifically the sparsity pattern of the divergence operator when using closed nodes is a strict subset of that when using half-closed nodes, whilst for the sparsity pattern of the LDG Laplacian is identical when using either closed or half-closed nodes. 

We also note that the sparsity pattern of the divergence/gradient operators for any choice of nodes to be strictly contained within the sparsity pattern of the LDG Laplacian. This is the as the LDG Laplacian is formed via a matrix product involving these operators. Consequently any linear solver technique dependent solely on the sparsity pattern of the Laplacian can also be applied to the gradient and divergence operators. As such we focus our discussion solely on solvers which leverage the sparsity structure of the Laplacian.

\subsection{Static condensation}
The first solver technique we consider is static condensation or Guyan reduction, popular in Finite Element methods (FEM). The idea behind static condensation is to via the Schur complement eliminate a number of unknowns from the linear system such that a smaller resultant system is obtained. For our construction here we focus on a strategy which uses the switch function to find an optimal elimination pattern within each element.

The first step of static condensation is to partition the set of nodes in a mesh into dependent and independent sets, where the dependent nodes are to be eliminated from the system. Given a linear system
\begin{equation}
    A u = f
\end{equation}
where $A$ is for instance the divergence or LDG Laplace operator described in Sect. \ref{sect:sparsity}, the nodes are reordered such that the independent nodes are ordered first and the linear system partitioned into the form
\begin{equation}
    \begin{bmatrix}
        A_{II} & A_{ID} \\
        A_{DI} & A_{DD}
    \end{bmatrix}
    \begin{bmatrix}
        u_I \\ u_D
    \end{bmatrix}
    = 
    \begin{bmatrix}
        f_I \\ f_D
    \end{bmatrix}
\end{equation}
where $u_I, u_D$ correspond to the values of a vector $u$ for the independent and dependent nodes respectively. Assuming now that the block corresponding only to the dependent nodes $A_{DD}$ can be easily inverted, a new reduced system can be obtained for only the independent nodes as
\begin{equation}
    \underbrace{\bigg( A_{II} - A_{ID} A_{DD}^{-1} A_{DI} \bigg)}_{\tilde{A}} u_I = \underbrace{f_I - A_{ID} A_{DD}^{-1} f_D}_{\tilde{f}}
\end{equation}
and we say that the dependent nodes have been eliminated, or condensed from the system. Once the solution for the independent nodes $u_I$ has been found, the solution $u_D$ on the dependent nodes can be easily found as the solution to
\begin{equation}
    A_{DD}u_D = f_D - A_{DI} u_I
\end{equation}

In particular if $A_{DD}$ is block diagonal, the Schur complement system is extremely easy both to form and to solve. Furthermore in this case the block sparsity pattern of $\tilde{A}$ will in fact remain to be the same as that of $A_{II}$. The goal therefore is to within each element find the largest possible set of dependent nodes that are eliminated such that $A_{DD}$ remains block diagonal.

To find this set of dependent nodes, we examine the sparsity patterns of the divergence and Laplace operators. If we consider two neighbouring elements $K_n \sim K_m$, as per the discussion in Sect. \ref{sect:sparsity} nodes across elements not on a boundary where the switch $S_n^m > 0$ do not communicate across elements. This implies simply that if we take the set of dependent nodes in each element $K_n$ to be the ones not on boundaries where the switch $S_n^m > 0$, the resulting submatrix $A_{DD}$ will be block diagonal.

For quadrilateral elements of polynomial degree $p$ then each element can be partitioned into $p^d$ dependent nodes to be eliminated and $(p+1)^d - p^d$ independent nodes in the reduced system. This is as in our methodology with Gauss-Radau nodes in each element exactly half the element boundaries will have boundary nodes placed on them. For simplex elements the number of dependent nodes that can be eliminated depends on the switch function in an element $K_n$, as unlike in quadrilateral case, the number of boundaries where the switch $S_n^m > 0$ differs from element to element.

\begin{figure}
    \centering
    \includegraphics[scale=0.4]{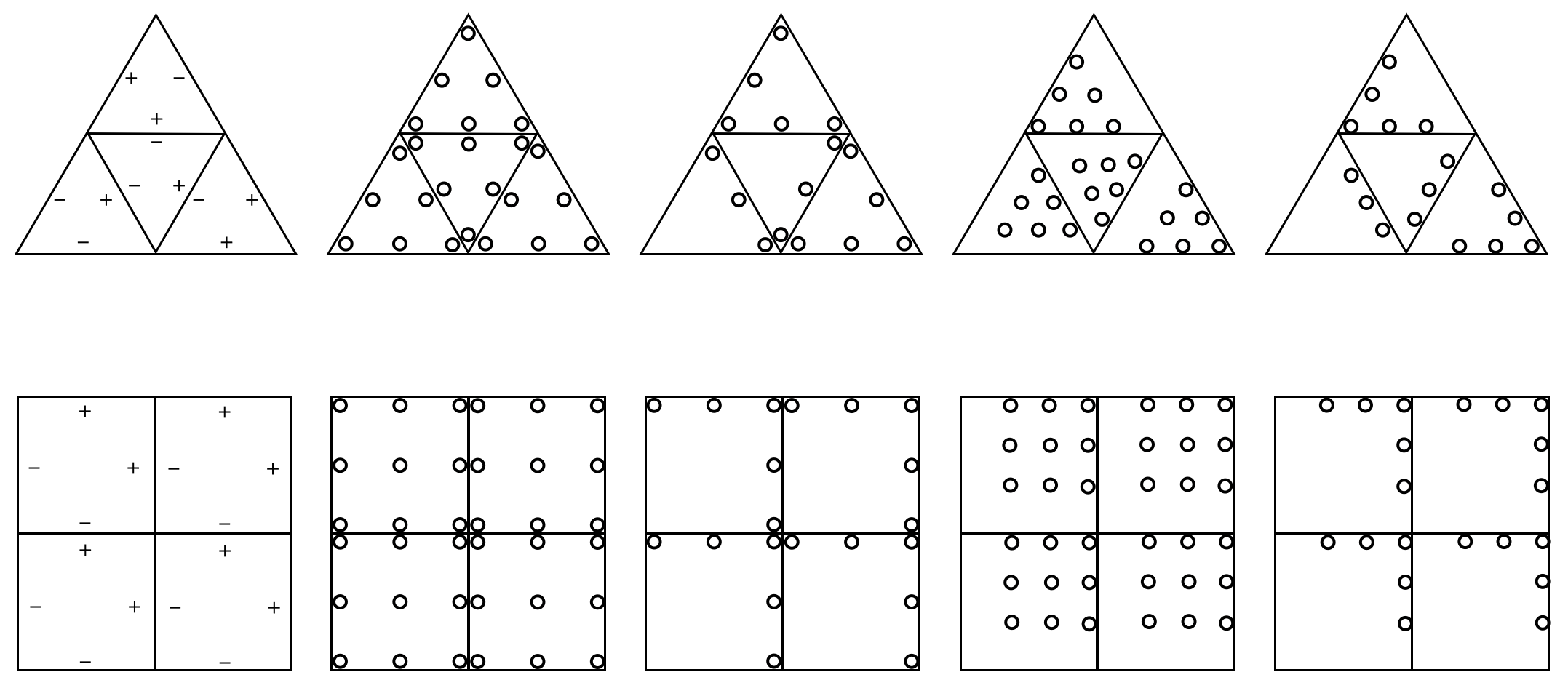}
    \caption{Static condensation elimination pattern in 2D for closed and half-closed nodes. In the leftmost column example switch functions on quadrilateral and simplex elements are shown. Column two shows $p=2$ closed nodes on the two meshes, and in column three nodes not on boundaries where the switch $S_n^m > 0$ are eliminated via static condensation. In columns four and five the same is shown for half-closed nodes.}
    \label{fig:elimination2d}
\end{figure}

Elimination according to this pattern guarantees that between two neighbouring elements $K_n \sim K_m$, only one of the elements will have dependent nodes placed on the boundary. This is because we require the switch function to always satisfy the property $S_n^m + S_m^n = 0$. This results in the reduced mesh having only $N_b$ boundary nodes left remaining on each boundary, which is similar to that of an eliminated degree $p+1$ Finite Element system of one degree higher. An example of the elimination procedure according to the switch function in 2D is illustrated in Fig. \ref{fig:elimination2d} on both a closed and half-closed mesh.

The idea of using static condensation for DG methods is in itself not new as has been used previously in the literature, for example in \cite{laskowski2020static,rueda2021static}. The technique is also related to the hybridized discontinuous Galerkin method \cite{cockburn2009unified}. However we are not aware of any elimination strategy that explicitly takes advantage of the sparsity pattern according to the switch function in this way, which allows for a greater number of dependent unknowns to be eliminated. Furthermore this strategy is easy to generalise to elements of arbitrary shape, since it only requires the definition of a valid switch function.

\subsubsection{1D example}
\begin{figure}
    \centering
    \includegraphics[scale=0.35]{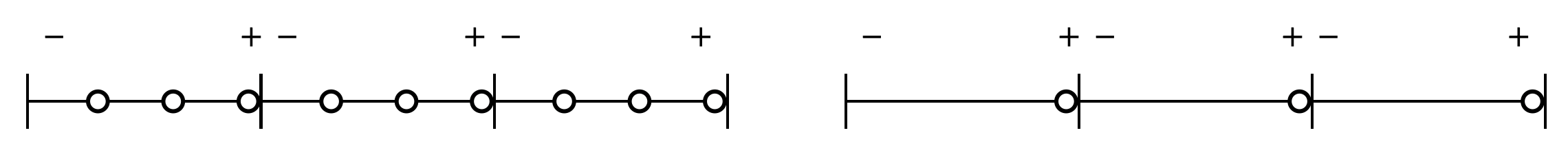}
    \caption{Static condensation for 1D example. Dependent nodes not on the boundary where the switch function $S_n^m > 0$ are eliminated to give a reduced operator on the remaining independent nodes.}
    \label{fig:elim1d}
\end{figure}

In the simple case of Poisson's problem on the 1D domain $\Omega = [0,1]$ with periodic boundary conditions, we can analyse the effect of static condensation more precisely. Without loss of generality we can freely divide the domain into $k$ elements of equal size with elements indexed from left to right and use the switch function given by the simple formula $S_n^m > 0$ if $n < m$ as shown in Fig. \ref{fig:elim1d}. For any given choice of closed or half-closed nodes according to this switch function, the nodes not on the right boundary of each element are designated as dependent nodes for the purposes of static condensation, and are eliminated from the LDG operator $L$ via a Schur complement to give a reduced system $\tilde{L}$.

Remarkably for any choice of closed or half-closed nodes, the resulting reduced operator $\tilde{L}$ has the same structure 
\begin{equation}
    \tilde{L}_{ii} = \frac{2}{h}, ~\tilde{L}_{i,i+1} = \tilde{L}_{i+1,i} = -\frac{1}{h}
\end{equation}
where $h = \frac{1}{k(p+1)}$. This is of course the well known centred Finite Difference approximation of the second derivative. This is a simple consequence of the following properties:
\begin{enumerate}
    \item The Schur complement of a symmetric positive matrix is symmetric positive semi-definite
    \item $\tilde{L}_{ii} = \tilde{L}_{jj}$ and $\tilde{L}_{i,i+1} = \tilde{L}_{j,j+1}$ for all $i,j$
    \item $\tilde{L}_{ij} \neq 0$ if and only if $j \in \{i-1,i,i+1\}$
    \item Rows of $\tilde{L}$ sum to $0$
\end{enumerate}
The first property is a well-known property of Schur complements, while the second is a consequence of symmetry. The third property can be deduced from the original sparsity pattern of $L$, as nodes in each element only communicate with those of its immediate neighbours in $L$. For the last property, we can simply note that the vector $E = (1,1,...,1)^T$ is an eigenvector of $L$ with eigenvalue $0$, and so
\begin{align*}
    &\begin{bmatrix}
        L_{II} & L_{ID} \\
        L_{DI} & L_{DD}
    \end{bmatrix}
    \begin{bmatrix}
        E_I \\ E_D
    \end{bmatrix}
    = 
    \begin{bmatrix}
        0 \\ 0
    \end{bmatrix}
    \\
    &\rightarrow~ L_{DD}^{-1}L_{DI} E_I = -E_D \\
    &\rightarrow~ L_{II}E_{I} - L_{ID}L_{DD}^{-1}L_{DI}E_I = 0
\end{align*}
which establishes the property.

\subsection{Block solvers}
A common type of solver used for DG problems are block based methods that leverage the block structure of the linear operators. These methods unlike the aforementioned static condensation can be applied on any DG method independent of the nodes used as they do not depend on the specific sparsity pattern. Popular examples of this include block Jacobi, block Gauss-Seidel and block-ILU solvers, which can either be used as a standalone method or as preconditioners or as smoothers for iterative methods. 

The cost of these methods is linked directly to the size of the blocks in the DG operators. This suggests a possible efficient linear solution strategy is to first apply static condensation to reduce the number of degrees of freedom in the linear system from $O(p^d)$ to $O(p^{d-1})$ before applying standard block based solvers to solve the reduced system. For many types of operators including the ones considered in this paper this is possible due to the fact that taking a Schur complement of a matrix preserves many properties of linear operators. In particular important properties such as strict diagonal dominance, symmetry and positive/negative definiteness of the eigenspectrum are all preserved under the Schur complement operation.

\subsection{Reduced block sizes}
We can analyse more quantitatively the block sizes of the statically condensed systems by counting the number of dependent degrees of freedom eliminated via the procedure compared to the total number of degrees of freedom in the system. As the number of dependent degrees of freedom are determined solely by the switch function, this number differs from simplex to quadrilateral meshes. In particular it is found that a higher proportion of the original nodes are eliminated in general for quadrilateral meshes compared to for simplex meshes though there is no difference in the asymptotic scaling of these proportions.

\subsubsection{Quadrilateral meshes}
For any quadrilateral mesh in $d$-dimensions for a fixed polynomial degree $p$, the number of unknowns in the DG problem is equal to $|\mathcal{T}_h|(p+1)^d$ where $|\mathcal{T}_h|$ denotes the number of elements. Following our elimination procedure, the number of unknowns that are eliminated is always equal to $|\mathcal{T}_h|p^d$. This is due to the switch function for all quadrilaterals satisfying the property of having opposite signs on opposite boundaries in an element, and that the elimination pattern is dependent only on the switch function. As such the proportion of unknowns that are eliminated in a quadrilateral mesh is always equal to $\frac{p^d}{(p+1)^d}$.

\subsubsection{Simplex meshes}
For a simplex mesh in $d$-dimensions for a fixed polynomial degree $p$, the number of unknowns in a DG formulation is equal to $|\mathcal{T}_h|{p+d \choose d}$ where $|\mathcal{T}_h|$ denotes the number of elements. However unlike with quadrilateral meshes, the number of unknowns eliminated with static condensation varies from element to element depending on its assigned switch function. Thus the exact number of nodes to be eliminated is problem dependent, specifically on the mesh and switch function used. However on a given mesh if assuming each boundary of an element to be equally likely to be assigned a switch of $\pm 1$, we are able to instead calculate an expected number of eliminated variables for a given polynomial degree $p$.

In 2D there are only two possibilities for a valid switch function assigned to a given element. In the first case two of its boundaries are assigned a switch of $+1$ with the remaining boundary assigned $-1$, whilst in the second case the reverse is true; these two configurations should be equally likely by assumption. Furthermore as shown in Fig. \ref{fig:elimination2d} it is easy to see that in each case, the number of dependent unknowns to be eliminated in each element $K_n$ is equal to ${p-b+d \choose d}$ where $b$ is the number of boundaries with a $S_n^m=+1$ switch assigned. This gives the expected proportion of degrees of freedom eliminated in 2D simplex mesh to be equal to $\frac{p^2}{(p+1)(p+2)}$.

The 3D case is slightly more complicated, in that there are three possible valid switch functions that may be assigned to a given element $K_n$. In the first case three of the four boundaries are assigned with a switch value of $S_n^m=+1$, in the second case two of the boundaries, and in the third and final case only one of the boundaries assigned this switch. By a simple counting argument it follows that each of these cases are assigned with probability $\frac{2}{7}, \frac{3}{7}, \frac{2}{7}$ respectively. Similar to the 2D case it can be seen that the number of dependent unknowns to be eliminated in each element is equal to ${p-b+d \choose d}$ where $b$ is the number of boundaries with a $S_n^m=+1$ switch assigned. Altogether then for an arbitrary 3D simplex mesh this gives the expected proportion of degrees of freedom eliminated by static condensation to be equal to $\frac{7p^3+5p}{7(p+1)(p+2)(p+3)}$.

\subsection{Static condensation + block solvers}
In addition to lowering the application cost of block based solvers, static condensation also improves the performance of the methods. Given a block preconditioner $P$ of a DG operator $A$, its efficacy can be measured by analysing the eigenspectrum of
\begin{equation} \label{eq:precond}
    I - P^{-1}A
\end{equation}
where the closer the eigenvalues to zero, the more effective the preconditioner. In the case where $P=A$, it can be seen that the eigenvalues of $I-P^{-1}A$ will all be exactly zero indicating a perfect preconditioner.

If static condensation is first performed on the matrix $A$ and the block preconditioner then applied onto the eliminated operator $\tilde{A}$, the operator whose eigenspectrum we need to analyse is instead
\begin{equation}
    \label{eq:condense_precond}
    I - 
    \begin{bmatrix}
        P_{\tilde{A}}^{-1} \tilde{A} & 0 \\
        C & I
    \end{bmatrix}
\end{equation}
where $P_{\tilde{A}}$ is the preconditioner constructed from $\tilde{A}$, and $C$ some arbitrary matrix whose value is unimportant for determining the spectrum of this operator. The derivation of this is provided in \ref{sect:appendix}. This implies that the degrees of freedom eliminated will all correspond to zero eigenvalues, implying improved solver performance.

\section{Numerical results}
In this section we perform numerical tests demonstrating the half-closed DG method on a few test problems. Examples demonstrating the accuracy using half-closed DG are shown, in addition to some focusing on the performance of the linear solvers described.

\subsection{Accuracy}
We apply here the method to quadrilateral meshes for the advection and Poisson equations. For these two equations, it has been previously proven in the literature that superconvergence of the solution is attained at the Gauss-Radau points when using Cartesian meshes with either Dirichlet or Neumann boundary conditions \cite{adjerid2002posteriori,liu2021superconvergence, yang2015analysis}. For a chosen polynomial degree $p$, it has been shown that order $p+2$ convergence is obtained at the Gauss-Radau points, as opposed to the standard $p+1$ convergence obtained by DG.

As the solution nodes on quadrilaterals are chosen to be the Gauss-Radau nodes, we perform a series of numerical tests examining the error obtained at the solution points. As the solution is expected to superconverge only at solution points, we do not use the continuous functional $L^2$ norm, but instead the discrete $L^2$ norm to analyse the error.

In each of the following examples, we compare the $L^2$ error attained at the nodes when using standard DG with closed Gauss-Lobatto nodes, compared to using half-closed DG with Gauss-Radau nodes. For each case, the examples were performed using three different choices of $p=1,2,3$, and the convergence behaviour for each polynomial degree $p$ is plotted separately.

\subsubsection{1D Advection equation}
We first consider the scalar advection equation in 1-dimension on the domain $\Omega = [0,1]$ with periodic boundary conditions on a uniform grid with $k$ elements:
\begin{equation}
    \frac{\partial}{\partial t} u + \alpha \cdot \nabla u = 0
\end{equation}
For this example, the initial condition was chosen to be $u(x) = \exp(100(x-\frac{1}{2})^2)$, and the final time was set to $T=1.0$, such that $u(0) = u(T)$. RK4 was chosen as the time integrator with a timestep of $10^{-6}$ employed to ensure convergence with respect to time. We examine the two cases of opposing signs of the velocity $\alpha$ separately, with the case of $\alpha < 1$ shown in Fig. \ref{fig:ex_downwind1d}, and of $\alpha > 1$ in Fig. \ref{fig:ex_upwind1d}.

\begin{figure}[h]
\centering
\pgfplotsset{every tick label/.append style={font=\huge}}
\begin{tikzpicture}[baseline=(current bounding box.center), 
    scale=0.5,transform shape]
      \begin{loglogaxis}[
            ymin=1e-3, ymax=1e0,
            xmin=6, xmax= 80,
            xtick={8, 16, 32, 64},
            xticklabels={$2^3$, $2^4$, $2^5$, $2^6$},
            ytick={1e0, 1e-1, 1e-2, 1e-3},
            yticklabels={$10^{0}$,$10^{-1}$,$10^{-2}$,$10^{-3}$,$10^{-4}$},
            grid = both,
            grid style = {line width=.1pt, draw=gray!15},
            major grid style = {line width=.2pt, draw=gray!50},
        ]
        \addplot[mark=none, color=black, mark size=3pt]
        table[ x=n, y=p1dg ]{conv1d_minus.txt};

        \addplot[mark=o, color=black, mark size=3pt]
        table[ x=n, y=p1hc ]{conv1d_minus.txt};

        \addplot[domain=30:60, samples=2] {80*2^(-2*log2(x))};
        \node at (axis cs:45, 7e-2) {\LARGE 2};

        \addplot[domain=30:60, samples=2] {400*2^(-3*log2(x))};
        \node at (axis cs:45, 2e-3) {\LARGE 3};

        \end{loglogaxis}
\end{tikzpicture}
\hspace{5mm}
\begin{tikzpicture}[baseline=(current bounding box.center), 
    scale=0.5,transform shape]
    \begin{loglogaxis}[
            ylabel style={ yshift=2ex },
            ymin=1e-6, ymax=1e0,
            xmin=6, xmax= 80,
            xtick={8, 16, 32, 64},
            xticklabels={$2^3$, $2^4$, $2^5$, $2^6$},
            ytick={1e0, 1e-2, 1e-4, 1e-6, 1e-8},
            yticklabels={$10^{0}$,$10^{-2}$,$10^{-4}$,$10^{-6}$,$10^{-8}$},
            grid = both,
            grid style = {line width=.1pt, draw=gray!15},
            major grid style = {line width=.2pt, draw=gray!50},
        ]
        \addplot[mark=none, color=black, mark size=3pt]
        table[ x=n, y=p2dg ]{conv1d_minus.txt};

        \addplot[mark=o, color=black, mark size=3pt]
        table[ x=n, y=p2hc ]{conv1d_minus.txt};

        \addplot[domain=30:60, samples=2] {200*2^(-3*log2(x))};
        \node at (axis cs:45, 7e-3) {\LARGE 3};

        \addplot[domain=30:60, samples=2] {50*2^(-4*log2(x))};
        \node at (axis cs:45, 3e-6) {\LARGE 4};

        \end{loglogaxis}
\end{tikzpicture}
\hspace{5mm}
\begin{tikzpicture}[baseline=(current bounding box.center), 
    scale=0.5,transform shape]
      \begin{loglogaxis}[
            ylabel style={ yshift=2ex },
            ymin=1e-8, ymax=1e0,
            xmin=6, xmax= 80,
            xtick={8, 16, 32, 64},
            xticklabels={$2^3$, $2^4$, $2^5$, $2^6$},
            ytick={1e0, 1e-2, 1e-4, 1e-6, 1e-8},
            yticklabels={$10^{0}$,$10^{-2}$,$10^{-4}$,$10^{-6}$,$10^{-8}$},
            grid = both,
            grid style = {line width=.1pt, draw=gray!15},
            major grid style = {line width=.2pt, draw=gray!50},
        ]
        \addplot[mark=none, color=black, mark size=3pt]
        table[ x=n, y=p3dg ]{conv1d_minus.txt};
        \label{plt:dg}

        \addplot[mark=o, color=black, mark size=3pt]
        table[ x=n, y=p3hc ]{conv1d_minus.txt};
        \label{plt:hc}

        \addplot[domain=30:60, samples=2] {700*2^(-4*log2(x))};
        \node at (axis cs:45, 1e-3) {\LARGE 4};

        \addplot[domain=30:60, samples=2] {350*2^(-5*log2(x))};
        \node at (axis cs:45, 6e-6) {\LARGE 5};

        \end{loglogaxis}
\end{tikzpicture}
\caption{Discrete $L^2$ error for 1D advection equation with negative velocity $\alpha=-1$. The $L^2$ error (y-axis) is plotted against the number of elements (x-axis) for polynomial orders $p=1,2,3$ from left to right. The error using DG with closed Gauss-Lobatto nodes is shown using \ref{plt:dg}, while for half-closed DG with Gauss-Radau nodes using \ref{plt:hc}.}
\label{fig:ex_downwind1d}
\end{figure}

\begin{figure}[h]
\centering
\pgfplotsset{every tick label/.append style={font=\huge}}
\begin{tikzpicture}[baseline=(current bounding box.center), 
    scale=0.5,transform shape]
      \begin{loglogaxis}[
            ymin=1e-3, ymax=1e0,
            xmin=6, xmax= 80,
            xtick={8, 16, 32, 64},
            xticklabels={$2^3$, $2^4$, $2^5$, $2^6$},
            ytick={1e0, 1e-1, 1e-2, 1e-3},
            yticklabels={$10^{0}$,$10^{-1}$,$10^{-2}$,$10^{-3}$,$10^{-4}$},
            grid = both,
            grid style = {line width=.1pt, draw=gray!15},
            major grid style = {line width=.2pt, draw=gray!50},
        ]
        \addplot[mark=none, color=black, mark size=3pt]
        table[ x=n, y=p1dg ]{conv1d_plus.txt};

        \addplot[mark=o, color=black, mark size=3pt]
        table[ x=n, y=p1hc ]{conv1d_plus.txt};

        \addplot[domain=30:60, samples=2] {80*2^(-2*log2(x))};
        \node at (axis cs:45, 7e-2) {\LARGE 2};
        
        \end{loglogaxis}
\end{tikzpicture}
\hspace{5mm}
\begin{tikzpicture}[baseline=(current bounding box.center), 
    scale=0.5,transform shape]
    \begin{loglogaxis}[
            ylabel style={ yshift=2ex },
            ymin=1e-6, ymax=1e0,
            xmin=6, xmax= 80,
            xtick={8, 16, 32, 64},
            xticklabels={$2^3$, $2^4$, $2^5$, $2^6$},
            ytick={1e0, 1e-2, 1e-4, 1e-6, 1e-8},
            yticklabels={$10^{0}$,$10^{-2}$,$10^{-4}$,$10^{-6}$,$10^{-8}$},
            grid = both,
            grid style = {line width=.1pt, draw=gray!15},
            major grid style = {line width=.2pt, draw=gray!50},
        ]
        \addplot[mark=none, color=black, mark size=3pt]
        table[ x=n, y=p2dg ]{conv1d_plus.txt};

        \addplot[mark=o, color=black, mark size=3pt]
        table[ x=n, y=p2hc ]{conv1d_plus.txt};

        \addplot[domain=30:60, samples=2] {200*2^(-3*log2(x))};
        \node at (axis cs:45, 7e-3) {\LARGE 3};
        
        \end{loglogaxis}
\end{tikzpicture}
\hspace{5mm}
\begin{tikzpicture}[baseline=(current bounding box.center), 
    scale=0.5,transform shape]
      \begin{loglogaxis}[
            ylabel style={ yshift=2ex },
            ymin=1e-8, ymax=1e0,
            xmin=6, xmax= 80,
            xtick={8, 16, 32, 64},
            xticklabels={$2^3$, $2^4$, $2^5$, $2^6$},
            ytick={1e0, 1e-2, 1e-4, 1e-6, 1e-8},
            yticklabels={$10^{0}$,$10^{-2}$,$10^{-4}$,$10^{-6}$,$10^{-8}$},
            grid = both,
            grid style = {line width=.1pt, draw=gray!15},
            major grid style = {line width=.2pt, draw=gray!50},
        ]
        \addplot[mark=none, color=black, mark size=3pt]
        table[ x=n, y=p3dg ]{conv1d_plus.txt};

        \addplot[mark=o, color=black, mark size=3pt]
        table[ x=n, y=p3hc ]{conv1d_plus.txt};

        \addplot[domain=30:60, samples=2] {700*2^(-4*log2(x))};
        \node at (axis cs:45, 1e-3) {\LARGE 4};
        
        \end{loglogaxis}
\end{tikzpicture}
\caption{Discrete $L^2$ error for 1D advection equation with positive velocity $\alpha=1$. The $L^2$ error (y-axis) is plotted against the number of elements (x-axis) for polynomial orders $p=1,2,3$ from left to right. The error using DG with closed Gauss-Lobatto nodes is shown using \ref{plt:dg}, while for half-closed DG with Gauss-Radau nodes using \ref{plt:hc}.}
\label{fig:ex_upwind1d}
\end{figure}

Interestingly the error obtained at solution nodes is superconvergent only if $\alpha < 1$. This in fact is consistent with the results presented in \cite{adjerid2002posteriori}, wherein superconvergence is only observed for hyperbolic equations if the velocity aligns with the Gauss-Radau nodes. Thus for a more general velocity in higher dimensions it is unclear as to whether there always exists a placement for the Gauss-Radau nodes where superconvergence is expected. To the author's knowledge there is also no analytical result establishing this as of yet, and we believe that it remains an open question.

\subsubsection{1D Poisson equation}
We next consider Poisson's equation in 1-dimension on the domain $\Omega = [0,1]$
\begin{equation}
    -\Delta u = f
\end{equation}
with Dirichlet boundary conditions on an uniform grid with $k$ elements. For the numerical test performed here, the solution is chosen as $u(x) = \exp( \sin(x) )$, with Dirichlet boundaries set according to this choice of solution and Dirichlet parameter $C_D=10/h$, where $h=\frac{1}{k}$ is the length an element in the domain. The errors in the discrete $L^2$ norm using standard DG with closed Gauss-Lobatto nodes and half-closed DG is shown in Fig. \ref{fig:ex_poisson1d}. We observe in this case that superconvergence is observed for all polynomial degrees $p$, consistent with the analytical findings in \cite{liu2021superconvergence}.

\begin{figure}[h]
\centering
\pgfplotsset{every tick label/.append style={font=\huge}}
\begin{tikzpicture}[baseline=(current bounding box.center), 
    scale=0.5,transform shape]
      \begin{loglogaxis}[
            ylabel style={ yshift=2ex },
            ymin=1e-5, ymax=1e0,
            xmin=6, xmax= 80,
            xtick={8,16,32,64},
            xticklabels={$2^3$,$2^4$,$2^5$,$2^6$},
            ytick={1e-0, 1e-2, 1e-4},
            yticklabels={$10^{0}$,$10^{-2}$,$10^{-4}$},
            grid = both,
            grid style = {line width=.1pt, draw=gray!15},
            major grid style = {line width=.2pt, draw=gray!50},
        ]
        \addplot[mark=none, color=black, mark size=3pt]
        table[ x=n, y=dg1 ]{poi1d.txt};

        \addplot[mark=o, color=black, mark size=3pt]
        table[ x=n, y=hc1 ]{poi1d.txt};

        \addplot[domain=30:60, samples=2] {50*2^(-2*log2(x))};
        \node at (axis cs:45, 7e-2) {\LARGE 2};

        \addplot[domain=30:60, samples=2] {10*2^(-3*log2(x))};
        \node at (axis cs:45, 3e-5) {\LARGE 3};
        
        \end{loglogaxis}
\end{tikzpicture}
\hspace{5mm}
\begin{tikzpicture}[baseline=(current bounding box.center), 
    scale=0.5,transform shape]
        \begin{loglogaxis}[
            ylabel style={ yshift=2ex },
            ymin=1e-8, ymax=1e-1,
            xmin=6, xmax= 80,
            xtick={8,16,32,64},
            xticklabels={$2^3$,$2^4$,$2^5$,$2^6$},
            ytick={1e-1, 1e-4, 1e-7},
            yticklabels={$10^{-1}$,$10^{-4}$,$10^{-7}$},
            grid = both,
            grid style = {line width=.1pt, draw=gray!15},
            major grid style = {line width=.2pt, draw=gray!50},
        ]
        \addplot[mark=none, color=black, mark size=3pt]
        table[ x=n, y=dg2 ]{poi1d.txt};

        \addplot[mark=o, color=black, mark size=3pt]
        table[ x=n, y=hc2 ]{poi1d.txt};

        \addplot[domain=30:60, samples=2] {80*2^(-3*log2(x))};
        \node at (axis cs:45, 3e-3) {\LARGE 3};

        \addplot[domain=30:60, samples=2] {1*2^(-4*log2(x))};
        \node at (axis cs:45, 7e-8) {\LARGE 4};
        
        \end{loglogaxis}
\end{tikzpicture}
\hspace{5mm}
\begin{tikzpicture}[baseline=(current bounding box.center), 
    scale=0.5,transform shape]
        \begin{loglogaxis}[
            ylabel style={ yshift=2ex },
            ymin=1e-10, ymax=1e-1,
            xmin=6, xmax= 80,
            xtick={8,16,32,64},
            xticklabels={$2^3$,$2^4$,$2^5$,$2^6$},
            ytick={1e-1, 1e-5, 1e-9},
            yticklabels={$10^{-1}$,$10^{-5}$,$10^{-9}$},
            grid = both,
            grid style = {line width=.1pt, draw=gray!15},
            major grid style = {line width=.2pt, draw=gray!50},
        ]
        \addplot[mark=none, color=black, mark size=3pt]
        table[ x=n, y=dg3 ]{poi1d.txt};

        \addplot[mark=o, color=black, mark size=3pt]
        table[ x=n, y=hc3 ]{poi1d.txt};

        \addplot[domain=30:60, samples=2] {200*2^(-4*log2(x))};
        \node at (axis cs:45, 3e-4) {\LARGE 4};

        \addplot[domain=30:60, samples=2] {20*2^(-5*log2(x))};
        \node at (axis cs:45, 5e-7) {\LARGE 5};
        
        \end{loglogaxis}
\end{tikzpicture}
\caption{Discrete $L^2$ error for 1D Poisson equation. The $L^2$ error (y-axis) is plotted against the number of elements (x-axis) for polynomial orders $p=1,2,3$ from left to right. The error using DG with closed Gauss-Lobatto nodes is shown using \ref{plt:dg}, while for half-closed DG with Gauss-Radau nodes using \ref{plt:hc}.}
\label{fig:ex_poisson1d}
\end{figure}

\subsubsection{2D Poisson equation}
\begin{figure}
    \centering
    \includegraphics[scale=0.25]{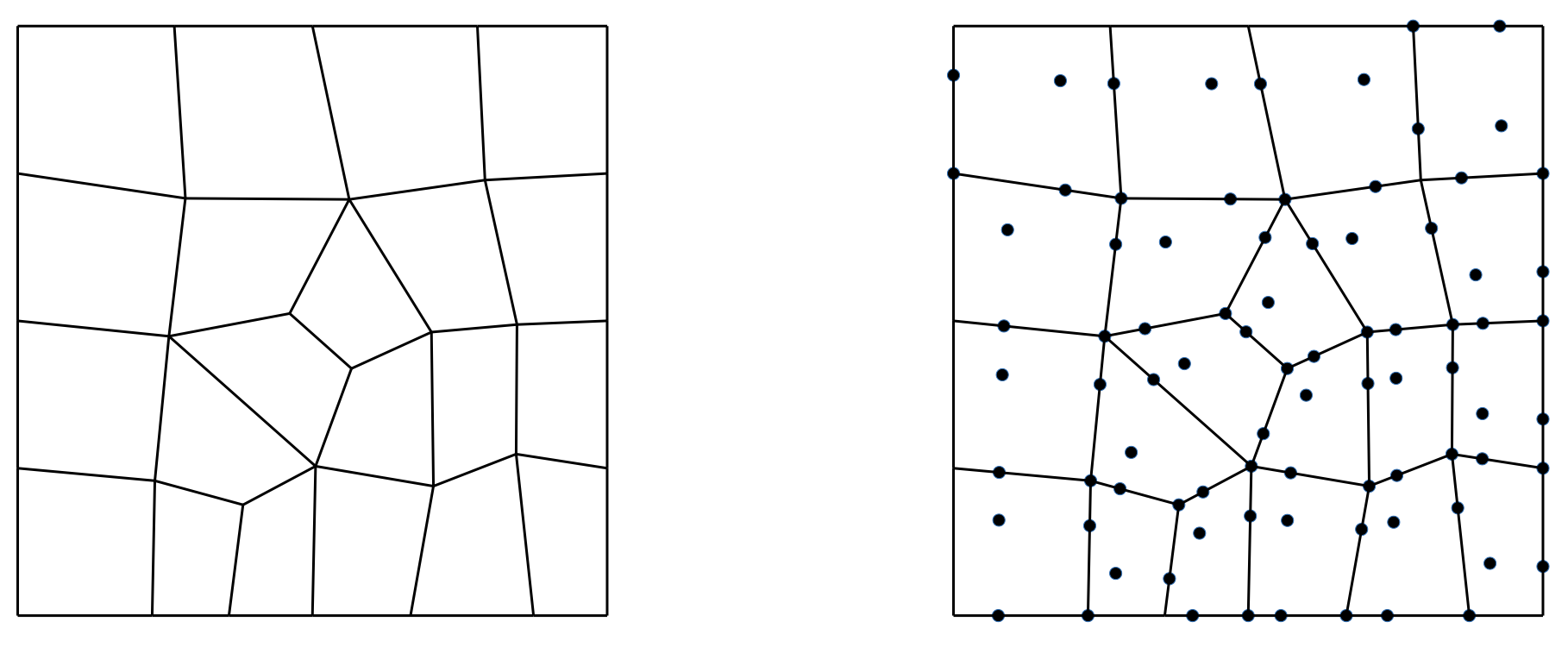}
    \caption{On the left, the mesh used for the 2D Poisson example at 0 refinements is shown. On the right the distribution of $p=1$ half-closed nodes on this mesh for this example is shown.}
    \label{fig:ex_mesh}
\end{figure}

As a final example we consider Poisson's equation in 2-dimension on the domain $\Omega = [0,1]^2$ with Dirichlet boundary conditions
\begin{equation}
    -\Delta u = f
\end{equation}
The mesh used for this test is shown in Fig. \ref{fig:ex_mesh}. For the numerical test performed here, the solution is chosen as $u(x,y) = \exp( \sin(x) \sin(y) )$, with Dirichlet boundaries set according to this choice of solution and Dirichlet parameter $C_D=10/h$, where $h$ is the mean length of elements on the domain boundary. The errors in the discrete $L^2$ using standard DG with closed Gauss-Lobatto nodes and half-closed DG with Gauss-Radau nodes is shown in Fig. \ref{fig:ex_poisson2d}.

\begin{figure}[ht]
\centering
\pgfplotsset{every tick label/.append style={font=\huge}}
\begin{tikzpicture}[baseline=(current bounding box.center), 
    scale=0.5,transform shape]
      \begin{semilogyaxis}[
            ylabel style={ yshift=2ex },
            ymin=1e-5, ymax=1e-0,
            xmin=-0.3, xmax= 3.3,
            xtick={0,1,2,3},
            xticklabels={0,1,2,3},
            ytick={1e-0, 1e-2, 1e-4},
            yticklabels={$10^{0}$,$10^{-2}$,$10^{-4}$},
            grid = both,
            grid style = {line width=.1pt, draw=gray!15},
            major grid style = {line width=.2pt, draw=gray!50},
        ]
        \addplot[mark=none, color=black, mark size=3pt]
        table[ x=nref, y=dg1 ]{poi2d.txt};

        \addplot[mark=o, color=black, mark size=3pt]
        table[ x=nref, y=hc1 ]{poi2d.txt};

        \addplot[domain=1.5:2.5, samples=2] {0.6*2^(-2*x)};
        \node at (axis cs:2, 9e-2) {\LARGE 2};

        \addplot[domain=1.5:2.5, samples=2] {0.015*2^(-3*x)};
        \node at (axis cs:2, 9e-5) {\LARGE 3};
        
        \end{semilogyaxis}
\end{tikzpicture}
\hspace{5mm}
\begin{tikzpicture}[baseline=(current bounding box.center), 
    scale=0.5,transform shape]
    \begin{semilogyaxis}[
            ylabel style={ yshift=2ex },
            ymin=1e-7, ymax=1e-0,
            xmin=-0.3, xmax= 3.3,
            xtick={0,1,2,3},
            xticklabels={0,1,2,3},
            ytick={1e0, 1e-2, 1e-4, 1e-6, 1e-8},
            yticklabels={$10^{0}$,$10^{-2}$,$10^{-4}$,$10^{-6}$,$10^{-8}$},
            grid = both,
            grid style = {line width=.1pt, draw=gray!15},
            major grid style = {line width=.2pt, draw=gray!50},
        ]
        \addplot[mark=none, color=black, mark size=3pt]
        table[ x=nref, y=dg2 ]{poi2d.txt};

        \addplot[mark=o, color=black, mark size=3pt]
        table[ x=nref, y=hc2 ]{poi2d.txt};

        \addplot[domain=1.5:2.5, samples=2] {0.2*2^(-3*x)};
        \node at (axis cs:2, 1e-2) {\LARGE 3};

        \addplot[domain=1.5:2.5, samples=2] {0.003*2^(-4*x)};
        \node at (axis cs:2, 3e-6) {\LARGE 4};
        
        \end{semilogyaxis}
\end{tikzpicture}
\hspace{5mm}
\begin{tikzpicture}[baseline=(current bounding box.center), 
    scale=0.5,transform shape]
      \begin{semilogyaxis}[
            ylabel style={ yshift=2ex },
            ymin=1e-8, ymax=1e-1,
            xmin=-0.3, xmax= 3.3,
            xtick={0,1,2,3},
            xticklabels={0,1,2,3},
            ytick={1e-1, 1e-3, 1e-5, 1e-7},
            yticklabels={$10^{-1}$,$10^{-3}$,$10^{-5}$,$10^{-7}$},
            grid = both,
            grid style = {line width=.1pt, draw=gray!15},
            major grid style = {line width=.2pt, draw=gray!50},
        ]
        \addplot[mark=none, color=black, mark size=3pt]
        table[ x=nref, y=dg3 ]{poi2d.txt};

        \addplot[mark=o, color=black, mark size=3pt]
        table[ x=nref, y=hc3 ]{poi2d.txt};

        \addplot[domain=1.5:2.5, samples=2] {0.05*2^(-4*x)};
        \node at (axis cs:2, 8e-4) {\LARGE 4};

        \addplot[domain=1.5:2.5, samples=2] {0.0005*2^(-5*x)};
        \node at (axis cs:2, 1e-7) {\LARGE 5};
        
        \end{semilogyaxis}
\end{tikzpicture}
\caption{Discrete $L^2$ error for 2D Poisson equation. The $L^2$ error (y-axis) is plotted against the number of refinements of the base mesh (x-axis) for polynomial orders $p=1,2,3$ from left to right. The error using DG with closed Gauss-Lobatto nodes is shown using \ref{plt:dg}, while for half-closed DG with Gauss-Radau nodes using \ref{plt:hc}.}
\label{fig:ex_poisson2d}.
\end{figure}

While the analytical result in \cite{yang2015analysis} only proves superconvergence on Cartesian grids, we regardless observe similar superconvergent behaviour in the $L^2$ norm in this example. While it is unclear as to whether this holds in general for all meshes, we have found similar behaviour for numerous other meshes in our testing. It however remains an open question as to when this behaviour can be expected, and also if it holds for a wider range of boundary conditions.

\subsection{Performance}
We perform some numerical experiments to demonstrate the application of the linear solvers outlined in this manuscript. For these examples we only consider the LDG Laplacian operator, in both $2-$ and $3-$ dimensions, with Dirichlet boundaries applied on all domain boundaries.

\subsubsection{Static condensation}

To study the effect of static condensation on reducing the problem complexity, we curated a set of meshes in $2-$ and $3-$ dimensions of both simplices and quadrilaterals. Ten meshes were selected for each group, for example $2$D quadrilaterals or $3$D simplices, and the LDG Laplacian formed on each mesh. Static condensation was then applied onto these matrices and the ratio of remaining degrees of freedom and remaining non-zero entries before and after the Schur complement operation was calculated for each mesh. For each group,  the obtained ratios were averaged over the set of ten meshes in the group. This was performed for polynomial orders $p=1$ through to $p=7$, and the results plotted in Fig. \ref{fig:static_condensation}.

\begin{figure}[ht]
\centering
\pgfplotsset{every tick label/.append style={font=\huge}}
\begin{tikzpicture}[baseline=(current bounding box.center), 
    scale=0.5,transform shape]
      \begin{axis}[
            ymin=0.0, ymax=1,
            xmin=0.5, xmax= 7.5,
            ymajorgrids,
            xtick pos=left,
            ytick pos=left
        ]
        \addplot[mark=triangle*, color=black, mark size=3pt]
        table[ x=p, y=Tet ]{condense_dof.txt};

        \addplot[mark=square*, color=black, mark size=3pt]
        table[ x=p, y=Hex ]{condense_dof.txt};

        \addplot[mark=triangle, color=black, mark size=3pt]
        table[ x=p, y=Tri ]{condense_dof.txt};

        \addplot[mark=square, color=black, mark size=3pt]
        table[ x=p, y=Quad ]{condense_dof.txt};
        
        \end{axis}
\end{tikzpicture}
\hspace{15mm}
\begin{tikzpicture}[baseline=(current bounding box.center), 
    scale=0.5,transform shape]
      \begin{axis}[
            ylabel style={ yshift=2ex },
            ymin=0.0, ymax=1,
            xmin=0.5, xmax= 7.5,
            ymajorgrids,
            xtick pos=left,
            ytick pos=left
        ]
        \addplot[mark=triangle*, color=black, mark size=3pt]
        table[ x=p, y=Tet ]{condense_nnz.txt};
        \label{plt:tet}

        \addplot[mark=square*, color=black, mark size=3pt]
        table[ x=p, y=Hex ]{condense_nnz.txt};
        \label{plt:hex}

        \addplot[mark=triangle, color=black, mark size=3pt]
        table[ x=p, y=Tri ]{condense_nnz.txt};
        \label{plt:tri}

        \addplot[mark=square, color=black, mark size=3pt]
        table[ x=p, y=Quad ]{condense_nnz.txt};
        \label{plt:quad}
        
        \end{axis}
\end{tikzpicture}
\caption{The ratio of degrees of freedom before and after static condensation (left), and of ratio of non-zeros before and after static condensation (right) are shown here. The study was performed for polynomial orders $p=1$ to $p=7$. The ratios for 3D tetrahedra are denoted as \ref{plt:tet}, 3D hexahedra as \ref{plt:hex}, 2D triangles as \ref{plt:tri}, and 2D quadrilaterals as \ref{plt:quad}.}
\label{fig:static_condensation}
\end{figure}
As expected, for higher polynomial degrees $p$, a higher proportion of degrees of freedom and non-zero entries are eliminated from the LDG Laplacian via static condensation. In general static condensation is also more effective at reducing these ratios on quadrilateral meshes than on simplex meshes in both $2-$ and $3-$ dimensions, consistent with our analysis.

\subsubsection{Static condensation + block solvers}
We also study the effect of preconditioning using block based solvers before and after static condensation. To do this we analyse the spectra of the operators in Eq. (\ref{eq:precond}) and Eq. (\ref{eq:condense_precond}) using two popular choices for the preconditioner $P$. The first is that of block Jacobi where $P=D$ is taken to consist of only the diagonal blocks of the DG operator $A$. The second is block Gauss-Seidel where $P=L$ is taken to be the lower block diagonal of the DG operator $A$. The spectra comparison for these cases is shown in Fig. \ref{fig:condense_precond}. For both of these cases the mesh shown in Fig. \ref{fig:ex_mesh} at one refinement with polynomial degree $p=2$ was used.

\begin{figure}[ht]
\centering
\pgfplotsset{every tick label/.append style={font=\huge}}
\begin{tikzpicture}[baseline=(current bounding box.center), 
    scale=0.5,transform shape]
      \begin{axis}[
            ylabel = {$\lambda$},
            label style={font=\Large},
            ymin=-1, ymax=1,
            xmin=0.5, xmax= 756,
            ymajorgrids,
            xtick={},
            xticklabels={},
            xtick style={draw=none}
        ]
        \addplot[mark=none, color=black]
        table[ x=N, y=jacobiFull ]{elimSmooth.txt};

        \addplot[mark=none, color=black, line width=2pt]
        table[ x=N, y=jacobiElim ]{elimSmooth.txt};
        
        \end{axis}
\end{tikzpicture}
\hspace{15mm}
\begin{tikzpicture}[baseline=(current bounding box.center), 
    scale=0.5,transform shape]
      \begin{axis}[
            ylabel = {$\lambda$},
            label style={font=\Large},
            ymin=-1, ymax=1,
            xmin=0.5, xmax= 756,
            ymajorgrids,
            xtick={},
            xticklabels={},
            xtick style={draw=none}
        ]
        \addplot[mark=none, color=black]
        table[ x=N, y=gsFull ]{elimSmooth.txt};

        \addplot[mark=none, color=black, line width=2pt]
        table[ x=N, y=gsElim ]{elimSmooth.txt};
        
        \end{axis}
\end{tikzpicture}
\caption{The spectra of the preconditioned system before and after static condensation using block Jacobi on the left and using block Gauss-Seidel on the right is shown. The thin line shows the spectra if the preconditioner is applied to the LDG Laplacian, while the bold line shows the spectra when the preconditioner is applied to the LDG Laplacian post static condensation.}
\label{fig:condense_precond}
\end{figure}

We clearly see that the number of zero eigenvalues in each case is greatly increased if preconditioning is performed post static condensation, as opposed to before. This implies for preconditioning to be more effective on the eliminated system than on the original as expected. In addition we emphasise that the preconditioner is cheaper to apply on the eliminated system than on the original system, which should also lower the overall cost of solving the linear system.

\section{Conclusion}
We have in this work introduced the concept of half-closed Discontinuous Galerkin discretisations, in which half-closed nodes are used to construct the nodal DG basis. For quadrilateral meshes we focus on the use of the Gauss-Radau nodes, which for straight sided elements gives a diagonal mass matrix. This was used to efficiently assemble other DG operators such as the divergence and LDG Laplace operators, as we showed that these can be cheaply assembled via left matrix multiplication with the mass matrix. However for simplex elements, if polynomial basis functions are desired, there do not exist any known nodes as of yet for arbitrary polynomial degrees such that a diagonal mass matrix is obtained.

It was found that using half-closed nodes does not increase the number of non-zero entries in first derivative DG operators significantly when compared to using closed nodes, and in the case of the Laplace operator to not increase it at all. Furthermore on quadrilateral meshes, half-closed nodes placed at the Gauss-Radau nodes are able to take advantage of previously established superconvergence results, such that a more accurate solution is obtained at the solution nodes.

As the sparsity pattern of the LDG Laplacian is identical whether using half-closed or closed nodes, we also discussed to use of some linear solver methods that can be applied when using either type of nodes. The methods we focused on were static condensation popular in Finite Element methods, and block preconditioners such as block Jacobi and block Gauss-Seidel popular in DG methods. We demonstrated that both these could be easily applied on both half-closed and closed DG discretisations, and could furthermore be used effectively in unison.

For future work we plan to apply the method on a range of more complicated test problems to measure more quantitatively how properties described in this work may translate to practical speedups for large problems. Another major avenue of interest is in exploring whether the extra degree of quadrature precision attained by the Gauss-Radau nodes in quadrilateral meshes provides any advantages in avoiding aliasing issues for nonlinear problems. Finally there is also the question of other choices of basis functions and nodes on simplices, such that some of the properties from using Gauss-Radau nodes on quadrilaterals might transfer over.

\section*{Acknowledgments}
The authors are grateful to Will Pazner for discussions which helped provide the inspiration for this work. This work was supported in part by the Director, Office of Science, Office of Advanced Scientific Computing Research, U.S. Department of Energy under Contract No. DE-AC02-05CH11231, and in part by the National Science Foundation under Grant DMS-2309596.

\newpage
\bibliographystyle{plain}
\bibliography{references}

\newpage
\appendix

\section{Preconditioning post static condensation}
\label{sect:appendix}
Given a preconditioner $P$ of a matrix $A$, an iterative method to solve $Ax=b$ can be constructed as
\begin{equation}
    x^{n+1} = (I - P^{-1} A) x^n + P^{-1}b
\end{equation}
where $x^n$ is the guess for the solution $x$ at the $n$-th iteration. For static condensation as before we write the system in terms of independent and dependent nodes
\begin{equation}
    \begin{bmatrix}
        A_{II} & A_{ID} \\
        A_{DI} & A_{DD}
    \end{bmatrix}
    \begin{bmatrix}
        u_I \\ u_D
    \end{bmatrix}
    = 
    \begin{bmatrix}
        f_I \\ f_D
    \end{bmatrix}
\end{equation}
and produce the eliminated system
\begin{equation}
    \underbrace{(A_{II} - A_{ID} A_{DD}^{-1} A_{DI}) }_{\tilde{A}}x_I = \underbrace{f_I - A_{ID} A_{DD}^{-1} f_D}_{\tilde{f}}
\end{equation}
We can now construct a preconditioner $P_{\tilde{A}}$ for this new system for which the iterative procedure becomes
\begin{equation}
    x^{n+1}_I = (I - P^{-1}_{\tilde{A}} \tilde{A}) x^n_I + P^{-1}_{\tilde{A}} \tilde{f}
\end{equation}
where at any point of the iteration, we can also recover the corresponding $x_D$ via solving the equation
\begin{equation}
    A_{DD}x_D^n = f_D - A_{DI}x_I^n 
\end{equation}
Overall this implies that preconditioning on the eliminated system produces the iterative procedure
\begin{equation}
    \begin{bmatrix}
        x_I^{n+1} \\
        x_D^{n+1}
    \end{bmatrix}
    =
    \begin{bmatrix}
        I-P^{-1}_{\tilde{A}}\tilde{A} & 0 \\
        A_{DI}(I-P^{-1}_{\tilde{A}})\tilde{A} & 0
    \end{bmatrix}
    \begin{bmatrix}
        x_I^{n} \\
        x_D^{n}
    \end{bmatrix}
    + 
    \begin{bmatrix}
        P^{-1}_{\tilde{A}} \tilde{f} \\
        A_{DD}^{-1}f_D
    \end{bmatrix}
\end{equation}
and so the convergence behaviour of the iterative method is determined by the spectrum of the matrix on the right hand side of the equation.

\end{document}